\documentclass[a4paper,12pt]{amsart}
\usepackage[foot]{amsaddr}
\usepackage{amssymb}
\usepackage{siunitx}
\usepackage{amsfonts}
\usepackage{booktabs}
\usepackage{amsmath,array}
\usepackage{caption}
\usepackage{subcaption}  
\usepackage{tabularx}
\usepackage{graphicx}
\usepackage[margin=2.7cm]{geometry}
\usepackage{bm}  
\usepackage{hyperref}

\begin{document}
	
	\title[Methods to enforce BCs in POD-Galerkin reduced order models]{A novel iterative penalty method to enforce boundary conditions in Finite Volume POD-Galerkin reduced order models for fluid dynamics problems}
	
	\author{S. Kelbij Star\textsuperscript{1,2*}}
	\address{\textsuperscript{1}SCK CEN, Institute for Advanced Nuclear Systems,Boeretang 200, 2400 Mol, Belgium.}
	\address{\textsuperscript{2}Ghent University, Department of Electromechanical, Systems and Metal Engineering, Sint-Pietersnieuwstraat 41, B-9000 Ghent, Belgium}
	\thanks{\textsuperscript{*}Corresponding Author.}
	\email{kelbij.star@sckcen.be}
	
	\author{Giovanni Stabile\textsuperscript{3}}
	\address{\textsuperscript{3}SISSA, International School for Advanced Studies, Mathematics Area, mathLab Trieste, Italy.}
	%\thanks{\textsuperscript{**}Second Corresponding Author.}
	\email{gstabile@sissa.it}
	
	%    author two information
	\author{Francesco Belloni\textsuperscript{1}}
	\email{francesco.belloni@sckcen.be}
	
	\author{Gianluigi Rozza\textsuperscript{3}}
	\email{grozza@sissa.it}
	
	\author{Joris Degroote\textsuperscript{2}}
	\email{Joris.Degroote@UGent.be}
	%\subjclass[2010]{78M34, 97N40, 35Q35}
	
	
	\date{}
	
	\dedicatory{}

	
	%%      ABSTRACT
	%%

%%%%% Begin Abstract %%%%%%%%%%%
\begin{abstract}
A Finite-Volume based POD-Galerkin reduced order model is developed for fluid dynamics problems where the (time-dependent) boundary conditions are controlled using two different boundary control strategies: the lifting function method, whose aim is to obtain homogeneous basis functions for the reduced basis space and the penalty method where the boundary conditions are enforced in the reduced order model using a penalty factor. The penalty method is improved by using an iterative solver for the determination of the penalty factor rather than tuning the factor with a sensitivity analysis or numerical experimentation. 

The boundary control methods are compared and tested for two cases: the classical lid driven cavity benchmark problem and a Y-junction flow case with two inlet channels and one outlet channel. The results show that the boundaries of the reduced order model can be controlled with the boundary control methods and the same order of accuracy is achieved for the velocity and pressure fields. Finally, the reduced order models are 270-308 times faster than the full order models for the lid driven cavity test case and 13-24 times for the Y-junction test case.
\end{abstract}

\maketitle

%%%% Start %%%%%%

\section{Introduction} 
Complex fluid dynamics problems are generally solved using discretization methods such as Finite Difference, Finite Element, Finite Volume (FV) or spectral element methods. However, it is usually not feasible to use these methods for applications that require to be solved almost in real time, such as on-the-spot decision making, (design) optimization or control~\cite{sartori2016reduced}. The high fidelity Computational Fluid Dynamics (CFD) tools, used for numerical simulations of the Navier--Stokes equations, are too computationally expensive for those purposes. This has motivated the development of reduced order modeling techniques. However, low degree-of-freedom models that are solely based on input-output data do not represent the physics of the underlying systems adequately and, moreover, may be sensitive to operating conditions~\cite{ravindran2000reduced}. 

Therefore, techniques, such as Reduced Basis (RB) methods, have been developed that retain the essential physics and dynamics of a high fidelity model that consists of discretized Partial Differential Equations (PDEs) describing the fluid problem~\cite{rozza2007reduced,veroy2003reduced}. The basic principle of these reduced order methods is to project the (parametrized) PDEs onto a low dimensional space, called the reduced basis space, in order to construct a physics-based model that is reduced in size and, therefore, in computational cost~\cite{hesthaven2016certified, quarteroni2015reduced, benner2015survey}. 

Fluid flows can be controlled in several ways. As an example, the system configuration can be manipulated by modifying the physical properties. However, in this work the focus is on controlling boundary conditions (BC) that are essential for defining flow problems. 

An example of a boundary control application from the nuclear field is the coupling of thermal-hydraulic system codes, i.e. transient simulations that are based on one-dimensional models of physical transport phenomena, with three-dimensional CFD codes~\cite{bandini2015assessment,toti2018coupled}. These type of system codes are, in general, based upon the solution of six balance equations for liquid and steam that are coupled with conduction heat transfer equations and that are supplemented by a suitable set of constitutive equations~\cite{petruzzi2008thermal}. 

One of the main purposes of this coupling is to speedup the CFD calculations by only including the region of interest in the CFD model and the rest of the domain in the much faster system code. However, the gain in computational time of such a coupled model is still limited by the CFD part. To overcome this burden, the system codes can be coupled with reduced order models (ROM) of the high fidelity CFD codes. For transient problems, time-dependent boundary conditions of the ROM are then to be controlled based on the BCs obtained from the systems codes.

For industrial applications, the Finite Volume discretization method is widely used by commercial software and open-source codes, as the method is robust~\cite{Eymard} and satisfies locally the conservation laws~\cite{versteeg2007, Fletcher}. 

By using a RB technique, the non-homogeneous BCs are, in general, no longer satisfied at the reduced order level. Furthermore, the BCs are not explicitly present in the ROM and therefore they cannot be controlled directly~\cite{lorenzi2016pod}. In literature~\cite{lorenzi2016pod,graham1999optimal1,kalashnikova2012efficient,Stabile2017CAIM}, different approaches to control the ROM BCs can be found of which two common approaches are extended and compared in this work: the lifting function method and the penalty method. The aim of the lifting function method~\cite{graham1999optimal1,Stabile2017CAIM} is to homogenize the BCs of the basis functions contained in the reduced subspace, while the penalty method~\cite{lorenzi2016pod,graham1999optimal1,kalashnikova2012efficient,Sirisup} weakly enforces the BCs in the ROM with a penalty factor. A disadvantage of the penalty method is that it relies on a penalty factor that has to be tuned with a sensitivity analysis or numerical experimentation~\cite{Sirisup}. Therefore, an iterative method is presented for tuning the penalty factor, which is, to the best of the authors' knowledge, introduced here for the first time in the context of Finite-Volume based POD-Galerkin reduced order methods. The novelty of this method is that a error tolerance for the enforced BC has to be set instead of an arbitrary value for the factor. Also the factor is determined automatically by iterating rather than manually via numerical experimentation.

The work is organized as follows: in Section~\ref{sec:FOM} the formulation of the full-order approximation of the PDEs is given and the methodology of the POD-based Galerkin projection is addressed in Section~\ref{sec:ROM}. In Section~\ref{sec:BCs} the two boundary control methods, the lifting function method and the iterative penalty method, are presented. In Section~\ref{sec:setup} the set-up of two numerical experiments, a lid driven cavity and a Y-junction test case, are given and the results are provided and discussed in Sections~\ref{sec:results} and~\ref{sec:discussion}, respectively. Finally, conclusions are drawn in Section~\ref{sec:conclusion} and an outlook for further developments is provided. 

\section{Full order model of the incompressible Navier--Stokes equations}\label{sec:FOM}
The fluid dynamics problem is physically described by the unsteady incompressible Navier--Stokes equations. In an Eulerian framework on a domain $\Omega$ $\subset$ $\mathbb{R}$$^d$ with $d$ = 2, 3 and boundary $\Gamma$ = ($\Gamma_{D_U} $ $\cup$ $\Gamma_{N_U} $) $\cap$ ($\Gamma_{D_p}$ $\cup$ $\Gamma_{N_p}$), the governing system of equations is given by 
\begin{align} \label{eq:FOM_mat}
\begin{cases}
\frac{\partial\boldsymbol{u}}{\partial t} + \nabla \cdot \left(\boldsymbol{u} \otimes \boldsymbol{u}\right)- \nabla \cdot (\nu \nabla \boldsymbol{u}) = - \nabla p + \boldsymbol{F}  &\mbox{in  } \Omega \times [0,T],  \\
\nabla \cdot \boldsymbol{u} = 0 &\mbox{in  } \Omega \times [0,T], \\
\boldsymbol{u}(\boldsymbol{x},0) = \boldsymbol{u}_0(\boldsymbol{x})&\mbox{in  } \Omega \times  \{0\}, \\
\boldsymbol{u} = \boldsymbol{f}(\boldsymbol{x},t) &\mbox{on  } \Gamma_{D_U}  \times [0,T], \\
\left(\nabla \boldsymbol{u} \right)\boldsymbol{n} = 0 &\mbox{on  } \Gamma_{N_U} \times [0,T],\\
\left(p\boldsymbol{I}\right)\boldsymbol{n} = 0 &\mbox{on  } \Gamma_{D_p} \times [0,T], \\
\left(\nabla p \right)\boldsymbol{n} = 0 &\mbox{on  } \Gamma_{N_p} \times [0,T],
\end{cases}
\end{align}
\noindent where $\boldsymbol{u}$ = $\boldsymbol{u}(\boldsymbol{x},t)$ represents the vectorial velocity field that is evaluated at $\boldsymbol{x} \in \Omega$ and $p = p(\boldsymbol{x},t)$ is the normalized scalar pressure field, which is divided by the constant fluid density $\rho$. $\nu$ is the kinematic viscosity and $\boldsymbol{F}$ is a body force term. For velocity, the (time-dependent) non-homogeneous Dirichlet boundary condition on $\Gamma_{D_U}$ is represented by  $\boldsymbol{f}(\boldsymbol{x},t)$ and $\boldsymbol{u}_0(\boldsymbol{x})$ denotes the initial condition for the velocity at time $t$ = 0 s. On $\Gamma_{N_U}$ a homogeneous Neumann boundary condition for velocity is applied and $\Gamma_{D_p}$ and $\Gamma_{N_p}$ are the Dirichlet and homogeneous Neumann boundary conditions for pressure. $\boldsymbol{n}$ denotes the outward pointing normal vector on the boundary and $T$ is the total simulation time.
The equations are presented here in a general format. The problem-specific (boundary) conditions are specified in Section~\ref{sec:setup}, in which the numerical experiments are presented.
\subsection{Pressure Poisson equation} \label{sec:PPE}
Standard Galerkin projection-based reduced order models are unreliable when applied to the non-linear unsteady Navier--Stokes equations~\cite{Lassila}. Furthermore, the ROMs need to be stabilized in order to produce satisfactory results for both the velocity and pressure fields~\cite{Sirisup,rozza2007stability,caiazzo2014numerical,Akhtar,bergmann2009enablers,noack2005need}. Two different stabilization techniques are compared in~\cite{Stabile2017CAF}; the supremizer enrichment of the velocity space in order to meet the inf-sup condition (SUP) and the exploitation of a pressure Poisson equation during the projection stage (PPE). The SUP-ROM performed about an order worse with respect to the PPE-ROM for what concerns the velocity field but better for what concerns the pressure field. This difference can be explained by the fact that within a supremizer stabilization technique, the POD velocity space is enriched by non-necessary (for the correct reproduction of the velocity field) supremizer modes. As the focus of this work is on controlling velocity boundary conditions, it is decided to use the PPE approach for stabilizing the ROM. Moreover, other approaches to simultaneously deal with velocity and pressure are the pressure stabilised Petrov--Galerkin methods~\cite{caiazzo2014numerical,baiges2014reduced,yano2014space} or assuming that velocity and pressure share the same temporal coefficients~\cite{lorenzi2016pod,bergmann2009enablers}.
For fluid problems that are solved numerically using a Finite Volume discretization technique~\cite{versteeg2007,moukalled2016finite}, often a Poisson Equation is solved for pressure as there is no dedicated equation for pressure in Equation~\ref{eq:FOM_mat}. The PPE is obtained by taking the divergence of the momentum equations and subsequently exploiting the divergence free constraint $\nabla \cdot \boldsymbol{u} = \boldsymbol{0}$. The resulting set of governing full order equations is then given by
\begin{align} \label{eq:PPE}
\begin{cases}
\frac{\partial\boldsymbol{u}}{\partial t} + \nabla \cdot \left(\boldsymbol{u} \otimes \boldsymbol{u}\right)- \nabla \cdot (\nu \nabla \boldsymbol{u}) = - \nabla p + \boldsymbol{F}  &\mbox{in  } \Omega \times [0,T], \\
\Delta p = - \nabla \cdot \left(\nabla \cdot \left(\boldsymbol{u} \otimes \boldsymbol{u} \right)\right) + \nabla \cdot \boldsymbol{F} &\mbox{in  } \Omega \times [0,T], \\
\boldsymbol{n} \cdot \nabla p = - \boldsymbol{n} \cdot  \left( \nu \nabla \times \nabla \times \boldsymbol{u} + \frac{\partial\boldsymbol{f}}{\partial t} \right) + \boldsymbol{n}  \cdot \boldsymbol{F} &\mbox{on  } \Gamma \times [0,T], \\ 
\end{cases}
\end{align}
In order to simplify the problem, no body force term, $\boldsymbol{F}$, is considered in this work. For more details on the derivation of the PPE the reader is referred to J.-G Liu et al.~\cite{liu2010stable}. These equations are discretized with the Finite Volume method and solved using a PIMPLE~\cite{ferziger2002computational} algorithm for the pressure-velocity coupling, which is a combination of SIMPLE~\cite{patankar1983calculation} and PISO~\cite{issa1986solution}.

\section{POD-Galerkin reduced order model  of the incompressible Navier--Stokes equations}\label{sec:ROM}
There exist several techniques in literature for creating a reduced basis space onto which the full order system (\ref{eq:FOM_mat}) is projected such as the Proper Orthogonal Decomposition (POD), the Proper Generalized Decomposition (PGD) and the Reduced Basis (RB) method with a greedy approach. For more details about the different methods the reader is referred to~\cite{rozza2007reduced,hesthaven2016certified,quarteroni2015reduced,chinesta2011short}. 
In this work, the Proper Orthogonal Decomposition method is used to create a reduced set of basis functions, or so-called modes, governing the essential dynamics of the full order model (FOM). For this, full order solutions are collected at certain time instances, the so-called snapshots. These snapshots do not necessarily have to be collected at every time step for which the full order solution is calculated.\\
Subsequently, it is assumed that the solution of the FOM can be expressed as a linear combination of spatial modes multiplied by time-dependent coefficients. The velocity snapshots $\boldsymbol{u}(\boldsymbol{x},t_n)$ and pressure snapshots $p(\boldsymbol{x},t_n)$ at time $t_n$ are approximated, respectively, by

\begin{equation}\label{eq:approx}
\boldsymbol{u}(\boldsymbol{x},t_n) \approx \boldsymbol{u}_r(\boldsymbol{x},t_n) = \sum\limits_{i=1}^{N_r^u} \boldsymbol{\varphi}_i(\boldsymbol{x})a_{i}(t_n), \hspace{0.5cm} p(\boldsymbol{x},t_n) \approx p_r(\boldsymbol{x},t_n) = \sum\limits_{i=1}^{N_r^p} \chi_i(\boldsymbol{x})b_{i}(t_n), 
\end{equation}
where $\boldsymbol{\varphi}_i$ and $\chi_i$ are the modes for velocity and pressure, respectively. \\$\boldsymbol{a}(t_n)$ = $\left[a_1(t_n), a_2(t_n), ..., a_2(t_n) \right]^T$ and $\boldsymbol{b}(t_n)$ = $ \left[b_1(t_n), b_2(t_n), ..., b_2(t_n) \right]^T$ are column vectors containing the corresponding time-dependent coefficients. $N_r^u$ is the number of velocity modes and $N_r^p$ is the number of pressure modes and thus it is assumed that velocity and pressure at reduced order level can be approximated with a different number of spatial modes. Furthermore, the modes are orthonormal to each other: ${\left( \boldsymbol{\varphi}_i,\boldsymbol{\varphi}_j\right)_{L^2(\Omega)}} = \delta_{ij}$, where $\delta$ is the Kronecker delta. The $L^2$-norm is preferred for discrete numerical schemes~\cite{Stabile2017CAF,busto2020pod} with ${\left( \cdot,\cdot\right)_{L^2(\Omega)}}$ the $L^2$-inner product of the fields over the  domain $\Omega$. \\
The optimal POD basis space for velocity, $E^{POD}_{u}$ = $\left[\boldsymbol{\varphi}_1,\boldsymbol{\varphi}_2, ... ,\boldsymbol{\varphi}_{N_r^u}\right]$ is then constructed by minimizing the difference between the snapshots and their orthogonal projection onto the basis for the $L^2$-norm~\cite{quarteroni2014reduced}. This gives the following minimization problem:
\begin{equation} \label{eq:min}
E^{POD}_{u} = \textrm{arg}\underset{\boldsymbol{\varphi}_1, ... ,\boldsymbol{\varphi}_{N_s^u}}{\textrm{min}} \frac{1}{N_s^u}\sum\limits_{n=1}^{N_s^u} \left\Vert \boldsymbol{u}(\boldsymbol{x},t_n) - \sum\limits_{i=1}^{N_r^u} \left( \boldsymbol{u}(\boldsymbol{x},t_n), \boldsymbol{\varphi_i} (\boldsymbol{x}) \right)_{L^2(\Omega)} \boldsymbol{\varphi}_i(\boldsymbol{x})\right\Vert_{L^2(\Omega)}^2,
\end{equation}
\noindent where $N_s^u$ is the number of collected velocity snapshots and $N_r^u$ (with $1\leq N_r^u$ $\leq$ $N_s^u$) denotes the dimension of the POD space $E^{POD}_{u}$. The POD modes are then obtained by solving the following eigenvalue problem on the snapshots~\cite{Stabile2017CAIM, Lassila,Stabile2017CAF,sirovich1987turbulence}:
\begin{equation} \label{eq:ev} 
\boldsymbol{C}\boldsymbol{Q}=\boldsymbol{Q}\boldsymbol{\lambda},
\end{equation}
\noindent where $C_{ij}$ = ${\left( \boldsymbol{u}(\boldsymbol{x},t_i),\boldsymbol{u}(\boldsymbol{x},t_j)\right)_{L^2(\Omega)}}$ for $i$,$j$ = 1, ..., $N_s^u$ is the correlation matrix, $\boldsymbol{Q}$ $\in \mathbb{R}^{N_s^u \times N_s^u}$ is a square matrix of eigenvectors and $\boldsymbol{\lambda}$ $\in \mathbb{R}^{N_s^u \times N_s^u}$ is a diagonal matrix containing the eigenvalues. The POD modes, $\boldsymbol{\varphi}_i$, are then constructed as follows
\begin{equation} \label{eq:POD}
\boldsymbol{\varphi}_i (\boldsymbol{x}) = \frac{1}{N_s^u\sqrt{\lambda_i}} \sum\limits_{n=1}^{N_s^u} \boldsymbol{u}(\boldsymbol{x},t_n) Q_{in}\text{\hspace{0.5cm} for \hspace{0.1cm}}i = 1,...,N_r^u,
\end{equation}
of which the most energetic (dominant) modes are selected. The procedure is the same for obtaining the pressure modes.\\
To obtain a reduced order model, the POD is combined with the Galerkin projection, for which the full order system of equations (Equation~\ref{eq:PPE}) is projected onto the reduced POD basis space. For more details about POD and Galerkin projection methods the reader is referred to~\cite{Stabile2017CAIM,Stabile2017CAF,georgaka2018parametric}. 
The following reduced system of momentum equations is then obtained 
\begin{equation}\label{eq:ROM}
%\begin{cases}
\boldsymbol{M_r} \dot{\boldsymbol{a}} +  \boldsymbol{C_r} (\boldsymbol{a}) \boldsymbol{a} - \nu \boldsymbol{A_r} \boldsymbol{a}  + \boldsymbol{B_r} \boldsymbol{b}  = 0,
%\end{cases}
\end{equation}
\noindent where the 'over-dot' indicates the time derivative and
\begin{equation}\label{eq:ROM_matrices}
\begin{split}
M_{r_{ij}}  = {\left( \boldsymbol{\varphi}_i,  \boldsymbol{\varphi}_j \right)_{L^{2}(\Omega)}}\text{\hspace{0.5cm} for \hspace{0.1cm}}i,j = 1,...,N_r^u, \\
A_{r_{ij}}  = {\left( \boldsymbol{\varphi}_i, \Delta \boldsymbol{\varphi}_j \right)_{L^{2}(\Omega)}}\text{\hspace{0.5cm} for \hspace{0.1cm}}i,j = 1,...,N_r^u, \\
B_{r_{ij}}  = {\left( \boldsymbol{\varphi}_i, \nabla \chi_j \right)_{L^{2}(\Omega)}}\text{\hspace{0.5cm} for \hspace{0.1cm}}i = 1,...,N_r^u \text{\hspace{0.1cm} and \hspace{0.1cm}}j = 1,...,N_r^p.
\end{split}
\end{equation}
These reduced matrices can be precomputed during an offline stage except for the non-linear term $\boldsymbol{C}_r$, which is given by
\begin{equation}\label{eq:C_matrix}
\begin{split}
C_{r_{ijk}}  = {\left( \boldsymbol{\varphi}_i, \nabla \cdot (\boldsymbol{\varphi}_j\otimes\boldsymbol{\varphi}_k) \right)_{L^{2}(\Omega)}}\text{\hspace{0.5cm} for \hspace{0.1cm}}i,j,k = 1,...,N_r^u, 
\end{split}
\end{equation}
This non-linear term is stored as a third order tensor~\cite{Stabile2017CAIM,quarteroni2007numerical} and the contribution of the convective term to the residual of Eq.~\ref{eq:ROM}, $R$, is evaluated at each iteration during the ROM simulations, or so-called online stage, as
\begin{equation}\label{eq:res}
R_i = \boldsymbol{a}^T C_{r_{i\bullet\bullet}}\boldsymbol{a}.
\end{equation}
The dimension of the tensor $\boldsymbol{C}_r$ (Equation~\ref{eq:C_matrix}) is growing with the cube of the number of modes used for the velocity space and therefore this approach may lead in some cases, especially when a large number of basis functions are employed, to high storage costs. Other approaches, such as EIM-DEIM~\cite{xiao2014non,barrault2004empirical} or Gappy-POD~\cite{carlberg2013gnat} may be more affordable~\cite{Stabile2017CAF}. 
As the pressure gradient term is present in the momentum equation the system is also coupled at reduced order level~\cite{Stabile2017CAIM}. The projection of the PPE leads to the following reduced system
\begin{equation}\label{eq:ROM2}
%\begin{cases}
\boldsymbol{D_r} \boldsymbol{b}  +  \boldsymbol{G_r}(\boldsymbol{a}) \boldsymbol{a} - \nu \boldsymbol{N_r} \boldsymbol{a} - \boldsymbol{T_r} \dot{\boldsymbol{a}}  = 0,
%\end{cases}
\end{equation}
\noindent where
\begin{equation}\label{eq:ROM_matrices2}
\begin{split}
D_{r_{ij}}  = {\left( \nabla \chi_i, \nabla \chi_j \right)_{L^{2}(\Omega)}} \text{\hspace{0.5cm} for \hspace{0.1cm}}i,j = 1,...,N_r^p, \\
G_{r_{ijk}}  = {\left( \nabla \chi_i, \nabla \cdot (\boldsymbol{\varphi}_j \otimes \boldsymbol{\varphi}_k ) \right)_{L^{2}(\Omega)}} \text{\hspace{0.5cm} for \hspace{0.1cm}}i = 1,...,N_r^p \text{\hspace{0.1cm} and \hspace{0.1cm}}j,k = 1,...,N_r^u, \\
N_{r_{ij}}  = {\left( \boldsymbol{n} \times \nabla \chi_i, \nabla \times \boldsymbol{\varphi}_j  \right)_{L^{2}(\Gamma)}}\text{\hspace{0.5cm} for \hspace{0.1cm}}i = 1,...,N_r^p \text{\hspace{0.1cm} and \hspace{0.1cm}}j = 1,...,N_r^u,\\
T_{r_{ij}}  = {\left( \chi_i, \boldsymbol{n} \cdot \boldsymbol{\varphi}_j  \right)_{L^{2}(\Gamma)}} \text{\hspace{0.5cm} for \hspace{0.1cm}}i = 1,...,N_r^p \text{\hspace{0.1cm} and \hspace{0.1cm}}j = 1,...,N_r^u,
\end{split}
\end{equation}
\noindent where the last two terms on the right hand side of Equation~\ref{eq:PPE} are projected on the boundary $\Gamma$.

Following the same strategy as in Equation~\ref{eq:res}, the non-linear term in Equation~\ref{eq:ROM2} is evaluated by storing the third order tensor $\boldsymbol{G_r}$. Equation~\ref{eq:ROM_matrices2} consists only of first order derivatives as integration by parts of the Laplacian term is used together with exploiting the pressure boundary condition after the PPE is projected onto the POD space spanned by the pressure modes. In that way, the numerical differentiation error can be reduced~\cite{Stabile2017CAF}. 

\subsection{Initial conditions} \label{sec:IC}
The initial conditions (IC) for the reduced system of Ordinary Differential Equations (Equations~\ref{eq:ROM} and~\ref{eq:ROM2}), are obtained by performing a Galerkin projection of the full order initial conditions onto the POD basis spaces as follows 
\begin{equation}\label{eq:ROM_IC}
\begin{split}
a_i(0) = \left(\boldsymbol{\varphi}_i,\boldsymbol{u}(0)\right)_{L^{2}(\Omega)}, \;\;\;\;\; b_i(0) = \left(\chi_i,p(0)\right)_{L^{2}(\Omega)},
\end{split}
\end{equation}
\noindent for velocity and pressure, respectively. It is important to note that the reduced system of equations are coupled and need to be solved iteratively. Moreover, the pressure is only defined up to an arbitrary constant, as in the FOM. Therefore, next to an initial condition for velocity, an initial guess for pressure is required for the system to converge more easily and to ensure the consistency between the FOM and the ROM~\cite{busto2020pod}.

In the online stage, the reduced system of equations (Equations~\ref{eq:ROM} and~\ref{eq:ROM2}) is solved for the velocity and pressure coefficients in the time period [$t_{1}$, $t_{\textrm{online}}$], where $t_{\textrm{online}}$ is the final simulation time.

\subsection{Relative error} \label{sec:error}
Three types of fields are considered: the full order fields, $X_{FOM}$, the projected fields, $X_r$, which are obtained by the $L^2$-projection of the snapshots onto the POD bases and the prediction fields obtained by solving the ROM, $X_{ROM}$. For every time instance, $t_n$, the basis projection error, $\|\hat{e}\|_{L^2}(\Omega)$, is given by
\begin{equation}\label{l2_prediction}
\|\hat{e}\|_{L^2(\Omega)}(t_n) = \frac{\|X_{FOM}(t_n)-X_{r}(t_n)\|_{L^{2}(\Omega)}}{\|X_{FOM}(t_n) \|_{L^{2}(\Omega)}},
\end{equation}
\noindent and the prediction error $\|e\|_{L^2}$, is determined by
\begin{equation}\label{l2_projection}
\|e\|_{L^2(\Omega)}(t_n) = \frac{\|X_{FOM}(t_n)-X_{ROM}(t_n)\|_{L^{2}(\Omega)}}{\|X_{FOM}(t_n) \|_{L^{2}(\Omega)}},
\end{equation}
\noindent where $X$ is either representing the velocity or pressure fields.

\section{Non-homogeneous (time-dependent) Dirichlet boundary conditions of the incompressible Navier--Stokes equations} \label{sec:BCs}
In a POD-based ROM, the non-homogeneous BCs are, in general, not satisfied by the ROM, as the basis functions, and in the same way their BCs, are a linear combination of the snapshots. Furthermore, the BCs are not explicitly present in the reduced system and therefore they cannot be controlled directly~\cite{lorenzi2016pod}. Two common approaches are presented in this section for handling the BCs: the lifting function and the penalty method~\cite{graham1999optimal1}. The aim of the lifting function method is to have homogeneous POD modes and to enforce the BCs by means of a properly chosen lifting function in the ROM. On the other side, the penalty method enforces the BCs in the ROM with a penalty factor. In this work only the velocity BCs are controlled with the two methods. 

\subsection{The lifting function method}\label{sec:control}
The lifting function method for the non-homogeneous boundary conditions is often used in the continuous Galerkin finite-element setting to reformulate a boundary control problem into a distributed one~\cite{chirco2019optimal,demkowicz2006computing,bornia2013distributed}. The method imposes the non-homogeneous (Dirichlet) conditions to the problem through lifting. This is done by subtracting the lifting function from the unknown variable in the original PDE problem, solving for the modified variable and adding the lifting function to the solution~\cite{ullmann2014pod}.
	
In a similar way, this method is used to impose non-homogeneous (Dirichlet) boundary conditions in reduced order models for which the lifted fields are projected onto the reduced bases spanned by the POD modes~\cite{fick2018stabilized}.

In this work, the velocity snapshots are made homogeneous by subtracting suitable lifting functions from all of them on which then the POD is performed. The result is a set of velocity modes that individually fulfill the homogeneous BCs as they are linear combinations of the modified velocity snapshots. The lifting functions, which fulfill the original non-homogeneous boundary conditions, are then added to a linear combination of POD basis functions. As a result, the non-homogeneous Dirichlet boundary conditions are included in the reduced basis space spanned by the POD modes and the lifting functions.

This lifting function method is also known as the "control function method" in literature~\cite{graham1999optimal1,Akhtar,lasiecka1984ritz} for PDE problems whose Dirichlet conditions can be parametrized with a single time-dependent coefficient~\cite{ullmann2014pod}. This is the type of problem that is presented in this work. The method is generalized in ~\cite{gunzburger2007reduced} for generic functions with multiple parameters at distinct boundary sections.

The functions to be chosen are system-specific and they have to satisfy the divergence free constraint in order to retain the divergence-free property of the snapshots~\cite{Stabile2017CAIM}. 

One way to generate a lifting function, $\tilde{\boldsymbol{\zeta}}_{c}(\boldsymbol{x})$ is by solving a problem as close as possible to the full order problem, where the boundary of interest is set to its value and everywhere else to a homogeneous BC. There are several other ways to compute the lifting function. For instance, the snapshot average can be used, although this does not always lead to a discretely divergence-free function. Alternatively, the solution of the stationary version of the considered problem can be computed~\cite{burkardt2006pod}. Two other common approaches are solving a non-homogeneous Stokes problem~\cite{fick2018stabilized,girault1999analysis,fonn2019fast} or solving a potential flow problem~\cite{eftang2010evaluation,HijaziStabileMolaRozza2020a}.

As one of the characteristics of the POD modes is that they are orthonormal, the lifting functions are normalized as follows
\begin{equation}\label{eq:norm_lift}
\boldsymbol{\zeta}_c(\boldsymbol{x}) = \frac{\tilde{\boldsymbol{\zeta}}_c(\boldsymbol{x})}{\|\tilde{\boldsymbol{\zeta}}_{c}(\boldsymbol{x})\|_{L^{2}(\Omega)}},
\end{equation}

before subtracting them from all snapshots and applying POD. The snapshots are then modified accordingly
\begin{equation}\label{eq:control}
\boldsymbol{u}^\prime(\boldsymbol{x},t)=\boldsymbol{u}(\boldsymbol{x},t)-\sum_{j=1}^{N_{BC}}\boldsymbol{\zeta}_{c_j}(\boldsymbol{x})u_{BC_j}(t) ,
\end{equation}
\noindent where $N_{BC}$ is the number of non-homogeneous BCs, $\boldsymbol{\zeta}_{c}(\boldsymbol{x})$ the normalized lifting functions and $u_{BC}$ is the normalized value of the corresponding Dirichlet boundary condition. 

The POD modes, $\boldsymbol{\varphi}^\prime_i$, that satisfy the homogeneous boundary conditions are obtained by solving an eigenvalue problem similar to Equation~\ref{eq:ev} on the homogenized snapshots $\boldsymbol{u}^\prime(\boldsymbol{x},t)$. The control functions are then added as additional modes to the reduced velocity basis
\begin{equation}\label{eq:approx_phi_lifted}
E^{\prime}_{u} = \left[\boldsymbol{\zeta}_{c_1}, ...,\boldsymbol{\zeta}_{c_{N_{BC}}}, \boldsymbol{\varphi}^\prime_1, ..., \boldsymbol{\varphi}^\prime_{N_r^u}  \right].
\end{equation}
Consequently, the velocity field at time $t_n$ is approximated by
\begin{equation}\label{eq:approx_temp}
\boldsymbol{u}_r(\boldsymbol{x},t_n)=\sum_{j=1}^{N_{BC}}\boldsymbol{\zeta}_{c_j}(\boldsymbol{x})u_{BC_j}(t_n) + \sum_{i=1}^{N_r^u}\boldsymbol{\varphi}^\prime_i(\boldsymbol{x})a_i(t_n),
\end{equation}
which satisfies the boundary conditions of the problem. $u_{BC}$ can be time-dependent. The Dirichlet boundary condition can be parametrized by assigning a new value to $u_{BC}$ in Equation~\ref{eq:approx_temp}. In other words, the lifting functions can be scaled by a factor.

For more details on the lifting function the reader is referred to~\cite{Stabile2017CAIM,georgaka2018parametric}. The overall algorithm for the lifting function method is given below. 

\clearpage

\begin{tabularx}{0.95\textwidth}{X}
	\toprule \textbf{Algorithm 1: lifting function method} \\\midrule
	\hspace{0.0cm}\textbf{OFFLINE PHASE}\\
	\hspace{0.0cm}\textbf{Solve full order model:}\\
	\hspace{0.0cm}(1) Generate snapshots over a time period [0, $T$] by solving the full order problem of Eq.~\ref{eq:FOM_mat};\\
	\hspace{0.0cm}\textbf{Obtain the lifting functions:}\\
	\hspace{0.0cm}(2) Generate the lifting functions by solving a flow problem: \\
	\hspace{0.6cm} \textbf{for} $i$ = 1 to $N_{BC}$ \textbf{do} \\
	\hspace{1.3cm} \textbf{for} $j$ = 1 to $N_{BC}$ \textbf{do} \\
	\hspace{2.0cm} \textbf{if} $i$ = $j$ \textbf{then}
	\\ \hspace{2.7cm}$\boldsymbol{u}$$\vert$$\Gamma_{D_j}$  = 1 \\
	\hspace{2.0cm} \textbf{else} \\
	\hspace{2.7cm}$\boldsymbol{u}$$\vert$$\Gamma_{D_j}$ = 0\\
	\hspace{1.3cm} \textbf{end for}\\
	\hspace{1.3cm} Solve a flow problem for $\tilde{\boldsymbol{\zeta}}_{c_i}$\\
	\hspace{0.6cm} \textbf{end} \textbf{for};\\
	\hspace{0.0cm}(3) Normalize the lifting functions to obtain $\boldsymbol{\zeta}_{c}$ as in Eq.~\ref{eq:norm_lift};\\
	\hspace{0.0cm}(4) Subtract the normalized lifting functions from the velocity snapshots as in Eq.~\ref{eq:control};\\
	\hspace{0.0cm}\textbf{Perform POD:}\\
	\hspace{0.0cm}(5) Retrieve the correlation matrix $\boldsymbol{C}$ from the solution snapshots;\\
	\hspace{0.0cm}(6) Solve the eigenvalue problem of Eq.~\ref{eq:ev} to obtain the POD modes using Eq.~\ref{eq:POD};\\
	\hspace{0.0cm}(7) Add the normalized lifting functions, $\boldsymbol{\zeta}_{c}$, as additional modes to the set of velocity POD modes $\boldsymbol{\varphi}$ according to Eq.~\ref{eq:approx_phi_lifted};\\
	\hspace{0.0cm}\textbf{Projection:}\\
	\hspace{0.0cm}(8) Project the discretized full order system onto the obtained reduced bases as done in Eq.~\ref{eq:ROM}-~\ref{eq:ROM_matrices2};\\
	\vspace{0.1cm}
	\hspace{0.0cm}\textbf{ONLINE PHASE}\\
	\hspace{0.0cm}\textbf{Solve reduced order model:}\\
	\hspace{0.0cm}(9) Project the initial fields for the parametrized BC onto the POD bases to get the initial condition/guesses for the ROM according to Eq.~\ref{eq:ROM_IC};\\
	\hspace{0.0cm}(10) Solve the reduced order problem of Eq.~\ref{eq:ROM} with the reduced Poisson equation, Eq.~\ref{eq:ROM2}, for pressure in the time period [$t_{1}$, $t_{\textrm{online}}$]; \\
	\hspace{0.0cm}(11) Reconstruct the full order fields from the obtained coefficients using Eq.~\ref{eq:approx_temp};
	\\\bottomrule
\end{tabularx}

\subsection{The iterative penalty method}
The penalty method was originally proposed in the context of finite element methods~\cite{lions1973non,babuvska1973finite}. The method transforms a strong non-homogeneous Dirichlet boundary condition into a weak Neumann boundary condition by the means of a small parameter whose inverse is called the penalty factor~\cite{placzek2008hybrid}. Thus, the method uses a penalty parameter to weakly impose the boundary conditions. In the POD-Galerkin reduced order modeling setting, the penalty method has been first introduced by Sirisup and Karniadakis~\cite{Sirisup} for the enforcement of boundary conditions.
For the penalty method, no modification of the snapshots is needed as the velocity Dirichlet BCs are directly enforced as constraints in the reduced system in the following way:
\begin{equation}\label{eq:pen_ROM_min}
\boldsymbol{M_r} \dot{\boldsymbol{a}}  - \nu \boldsymbol{B_r} \boldsymbol{a} + \boldsymbol{a}^T \boldsymbol{C_r} \boldsymbol{a} + \boldsymbol{K_r} \boldsymbol{b} +  \sum_{l=1}^{N_{BC}} \tau_l\left( \boldsymbol{P1}_l\boldsymbol{a} - u_{BC_l}(t) \boldsymbol{P2}_l \right) = 0 ,
\end{equation}
\noindent where $\tau$ is the penalty factor~\cite{Sirisup} and the additional terms with respect to Equation~\ref{eq:ROM} are projected on the boundary as follows
\begin{equation}\label{eq:D}
\begin{split}
P1_{lij} =  \left( \boldsymbol{\varphi}_i, \boldsymbol{\varphi}_j \right)_{L^{2}(\Gamma_l)}   \text{\hspace{0.5cm} for \hspace{0.1cm}}l = 1,...,N_{BC} \text{\hspace{0.1cm} and \hspace{0.1cm}} i,j = 1,...,N_r^u,\\
P2_{li}  = \left( \boldsymbol{\varphi}_i, \boldsymbol{\phi} \right)_{L^{2}(\Gamma_l)} \text{\hspace{0.5cm} for \hspace{0.1cm}}l = 1,...,N_{BC} \text{\hspace{0.1cm} and \hspace{0.1cm}} i = 1,...,N_r^u, 
\end{split}
\end{equation}
where $\boldsymbol{\phi}$ is a unit field. This minimization problem is formulated at reduced order level and, therefore, the penalty method does not depend on the full order snapshots. 

In order to have an asymptotically stable solution, the penalty factors $\tau$ should be larger than 0. If $\tau \rightarrow \infty$ the solution generally converges to a true optimal solution of the original unpenalized problem~\cite{hou1999numerical}. Nevertheless, a strong imposition would be approached and the ROM becomes ill-conditioned~\cite{lorenzi2016pod,epshteyn2007estimation}. Therefore, the penalty factor needs to be chosen above a threshold value for which the method is stable and converges~\cite{epshteyn2007estimation,dai2013analysis}. On the other hand, it is important to find a penalty factor as small as possible to obtain a numerical stable solution. This is usually done by numerical experimentation~\cite{lorenzi2016pod,graham1999optimal1,kalashnikova2012efficient, bizon2012reduced}. 

Several techniques exist in literature to optimize the numerical experimentation. Kelley~\cite[page~214]{kelley1962method} used a simple iteration scheme to optimize the trial-and-error process of the numerical experimentation. With this scheme the penalty value is adjusted each iteration by using the absolute value of the ratio between the constraint violation and a preassigned tolerance as a factor to increase or decrease the values at the end of each iteration. Basically, the idea is that the penalty factor obtained by the iteration scheme is optimal in the sense that it perturbs the original problem by a minimum for the given tolerance~\cite{huettelminimum}.

In this work the experimentation is optimized using a first-order iterative optimization scheme~\cite{leitmann1962optimization} to determine the factors that is based on the iteration scheme described in the previous paragraph. The penalty factors, $\tau$, are updated each iteration $k$, as follows
\begin{equation}\label{eq:tau}
\tau_l^{k+1}(t_n) = \tau_l^k (t_n) \frac{\left|r_l^k(t_n)\right|}{\epsilon} = \tau_l^k (t_n) \frac{\left|\tilde{u}_{BC_l}^k(t_n)-u_{BC_l}(t_n)\right|}{\epsilon} \text{\hspace{0.5cm} for \hspace{0.1cm}}l = 1,...,N_{BC},
\end{equation}
with $r^k(t_n)$ the residual between $\tilde{u}_{BC}^k$, the value of a certain boundary at the $k^{th}$ iteration, and $u_{BC}$, the enforced boundary condition, at an evaluated time $t_n$. $\tilde{u}_{BC}^k$ is obtained during the online phase by reconstructing the boundary. $\epsilon$ $>$ 0 is the given error tolerance for the residual which has to be set. There is no single approach that can be considered the best for choosing $\epsilon$, as the preferred tolerance depends on the problem and on both physical and geometrical parameters. The eigenvalue truncation error of the POD modes gives a good indication for the value of $\epsilon$. The penalty method is therefore no longer based on an arbitrary value for the penalty factor.

As long as $\left|\tilde{u}_{BC_l}^k(t_n)-u_{BC_l}(t_n)\right|$ $>$ $\epsilon$ the penalty factors grow every update and converge to the smallest penalty factors that satisfy the required tolerance. Thus, if the initial guess for the factor is below the minimum value for $\tau$ for which the boundary condition is enforced in the ROM, the factor is approached from below using this method. For a time-dependent problem it is not needed to determine a penalty factor for all time steps $N_t$. Often the factor determined after the first couple of time steps, $N_{\tau}$, can be used for the whole ROM solution.  

The step-by-step demonstration of the iterative function method is given below by Algorithm 2. 

\begin{tabularx}{0.95\textwidth}{X}
	\toprule \textbf{Algorithm 2: Iterative penalty method} \\\midrule
	\hspace{0.0cm}\textbf{OFFLINE PHASE}\\
	\hspace{0.0cm}\textbf{Solve full order model:}\\
	\hspace{0.0cm}(1) Generate snapshots over a time period [0, $T$] by solving the full order problem of Eq.~\ref{eq:FOM_mat};\\
	\hspace{0.0cm}\textbf{Perform POD:}\\
	\hspace{0.0cm}(2) Retrieve the correlation matrix $\boldsymbol{C}$ from the solutions;\\
	\hspace{0.0cm}(3) Solve the eigenvalue problem of Eq.~\ref{eq:ev} to obtain the POD modes using Eq.~\ref{eq:POD};\\
	\hspace{0.0cm}\textbf{Impose BCs with penalty method:}\\
	\hspace{0.0cm}(4) Project the modes on the reduced basis at the boundary of the domain to determine $\boldsymbol{P1}$ and $\boldsymbol{P2}$ for each non-homogeneous Dirichlet boundary condition as in Eq.~\ref{eq:D};\\
	\hspace{0.0cm}(5) Solve iteratively for the penalty factor using Eq.~\ref{eq:tau}:\\
	\hspace{0.6cm}\textbf{for} $i = 1$ to $N_{\tau}$ \textbf{do} \\
	\hspace{1.3cm}\textbf{while} $\left|\tilde{u}_{BC_l}^k(t_i)-u_{BC_l}(t_i)\right| > \epsilon$ \textbf{do}\\
	\hspace{2.0cm} $\tau_l^{k+1}(t_i)= \tau_l^k (t_i) \frac{\left|\tilde{u}_{BC_l}^k(t_i)-u_{BC_l}(t_i)\right|}{\epsilon}$ \text{\hspace{0.5cm}}\\
	\hspace{1.3cm}\textbf{end while}\\
	\hspace{0.6cm}\textbf{end for};\\
	\hspace{0.0cm}\textbf{Projection:}\\
	\hspace{0.0cm}(6) Project the discretized full order system onto the obtained reduced bases as done in Eq.~\ref{eq:ROM}-~\ref{eq:ROM_matrices2};\\
	\vspace{0.1cm}
	\hspace{0.0cm}\textbf{ONLINE PHASE}\\
	\hspace{0.0cm}\textbf{Solve reduced order model:}\\
	\hspace{0.0cm}(7) Project the initial fields for the parametrized BC onto the POD bases to get the initial condition/guesses for the ROM using Eq.~\ref{eq:ROM_IC};\\
	\hspace{0.0cm}(8) Solve the reduced order problem of Eq.~\ref{eq:ROM} with the reduced Poisson equation, Eq.~\ref{eq:ROM2}, for pressure in the time period [$t_{1}$, $t_{\textrm{online}}$]; \\
	\hspace{0.0cm}(9) Reconstruct the full order fields from the obtained coefficients using Eq.~\ref{eq:approx};
	\\\bottomrule
\end{tabularx}\\

It is important to note that the penalty factor can affect the number of iterations needed to solve the reduced system and therefore the convergence and cost of the reduced order model~\cite{nour1987note}.

\section{Numerical simulation tests} \label{sec:setup}
In this section the set-up of two cases are described for which the boundary control methods, the lifting function method and the iterative penalty method, are tested. The first one test case is the classical lid driven cavity benchmark problem and the second one is a Y-junction with two inlets and one outlet channel whose time-dependent inlet boundary conditions are controlled. 

\subsection{Lid-driven cavity flow problem}
\begin{figure}[h!]
	\centering
	\captionsetup{justification=centering}
	\includegraphics[width=8.0cm]{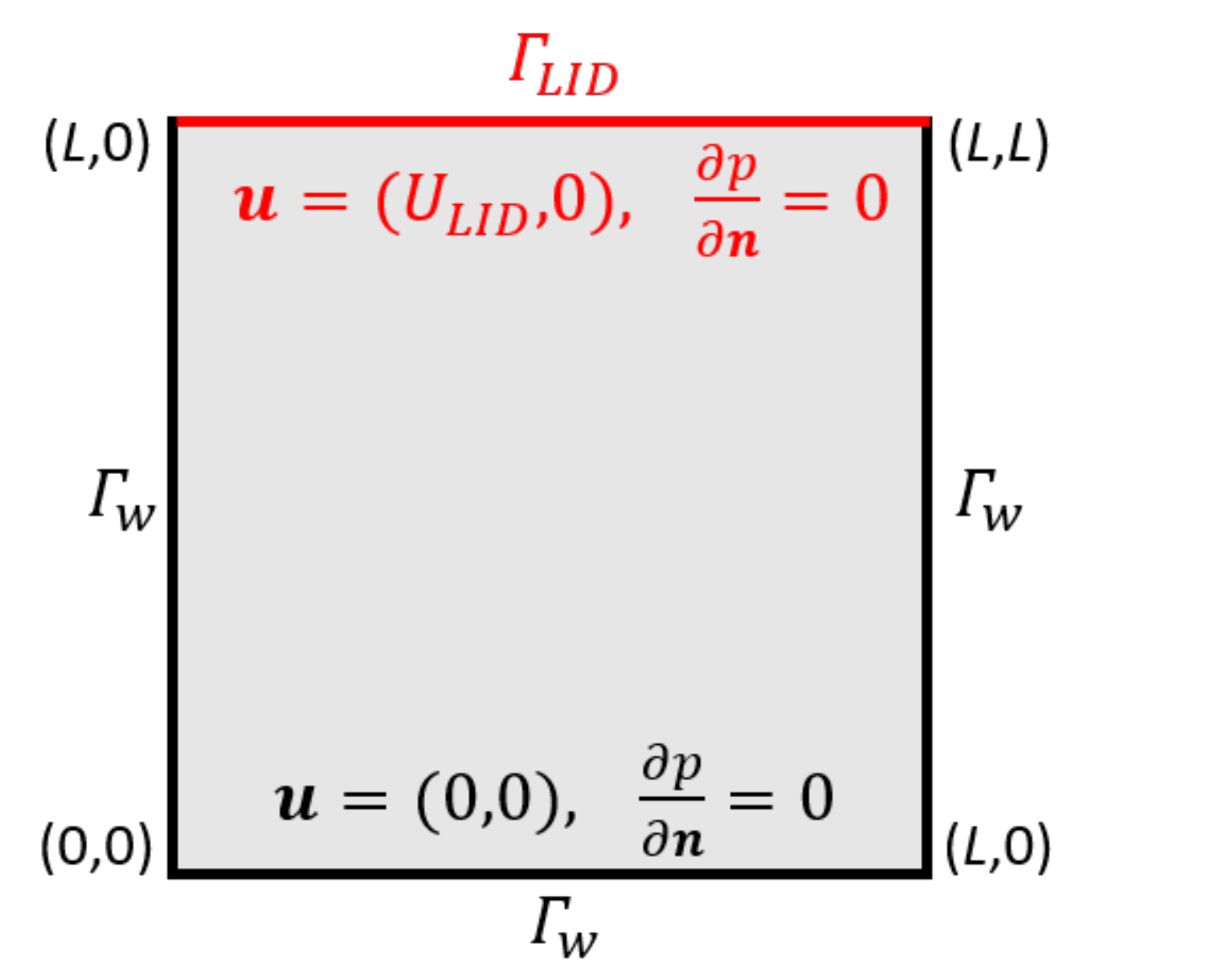}
	\caption{Sketch of the geometry of the 2D square cavity with moving top lid including boundary conditions.}
	\label{fig:LID_setup}
\end{figure}

The first test case consists of a lid driven cavity problem. The simulation is carried out on a two-dimensional square domain of length $L$ = 0.1 m on which a (200 $\times$ 200) structured mesh with quadrilateral cells is constructed. The boundary is subdivided into two different parts $\Gamma$ = $\Gamma_{LID}$ $\cup$ $\Gamma_w$ and the boundary conditions for velocity and pressure are set according to Figure~\ref{fig:LID_setup}. The pressure reference value is set to 0 m$^2$/s$^2$ at coordinate (0,0). At the top of the cavity a constant uniform and horizontal velocity equal to $\boldsymbol{u}$ = ($U_{LID}$,0) = (1,0) m/s is prescribed. A no slip BC is applied at the walls, $\Gamma_w$. The kinematic viscosity is equal to $\nu$ = 1 $\cdot$ \num{e-4} m$^2$/s and the corresponding Reynolds number is 1000, meaning that the flow is considered laminar. 

The unsteady full order equations are iteratively solved by the FV method with the $pimpleFoam$ solver of the open source C++ library OpenFOAM 6~\cite{Jasak}. The PIMPLE algorithm is used for the pressure-velocity coupling~\cite{ferziger2002computational}. For the full order simulations, the spatial discretization of all terms is performed with a central differencing scheme (linear). The temporal discretization is treated using a second order backward differencing scheme (BDF2). A constant time step of $\Delta t$ = 5 $\cdot$ \num{e-4} s has been applied and the total simulation time is 10 s. Snapshots of the velocity and pressure fields are collected every 0.01 s, resulting in a total of 1001 snapshots (including 1 for the initial condition). The initial condition field with $U_{LID}$ = 1 m/s is used as a lifting function. 

For this test case the same boundary conditions are applied in the ROM as in the FOM for which the snapshots are collected. The temporal discretization of the ROM is performed with a first order Newton's method.

POD, projection of the full order solution on the reduced subspace and the reduced order simulations are all carried out with ITHACA-FV, a C++ library based on the Finite Volume solver OpenFOAM. For more details on the ITHACA-FV code, the reader is referred to~\cite{Stabile2017CAIM,Stabile2017CAF,ITHACA}. 

\subsection{Y-junction flow problem}
\begin{figure}[h!]
	\centering
	\begin{subfigure}{.45\textwidth}
		\centering
		\includegraphics[width=1.\linewidth]{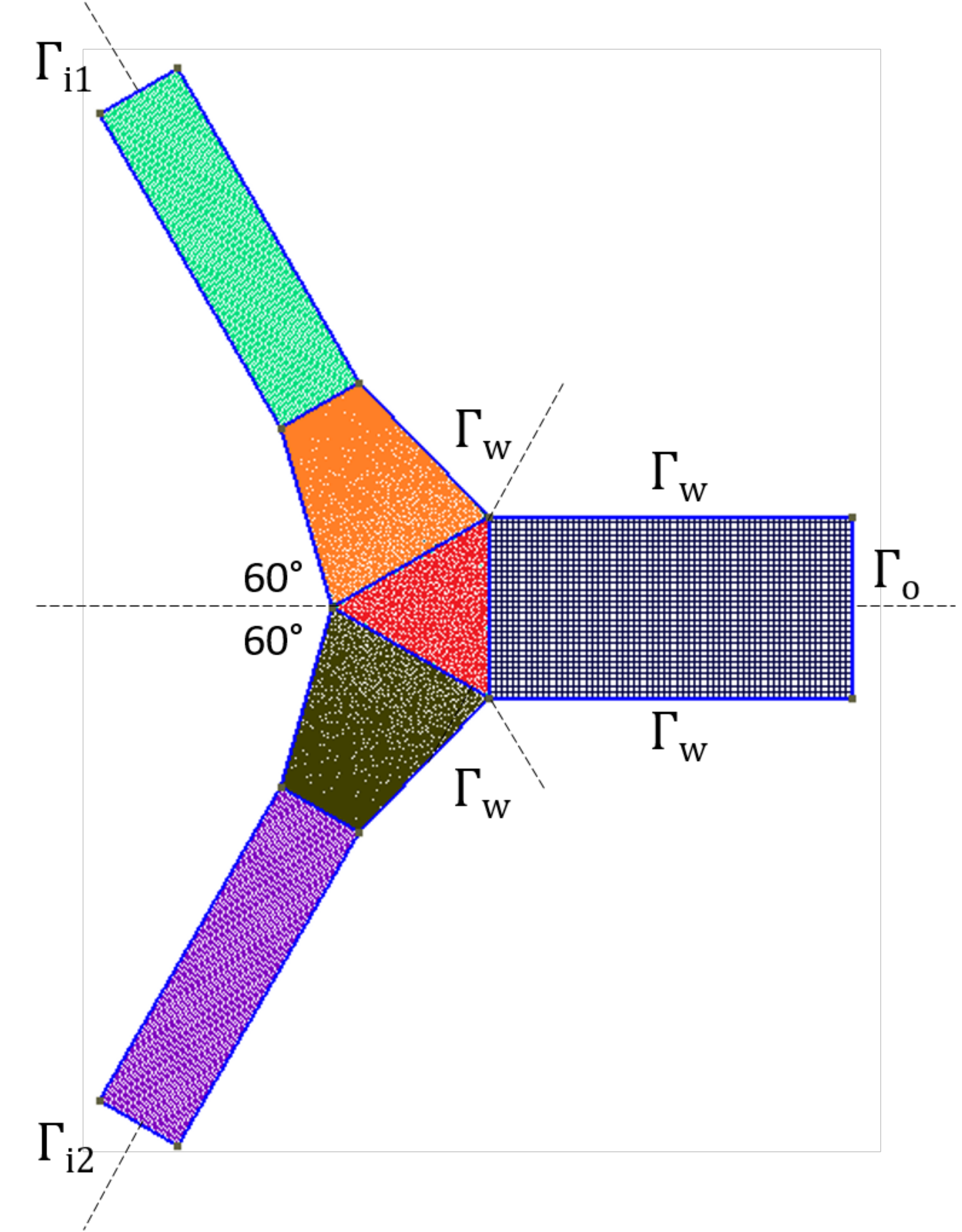}
	\end{subfigure}%
	\begin{subfigure}{.55\textwidth}
		\centering
		\includegraphics[width=1.\linewidth]{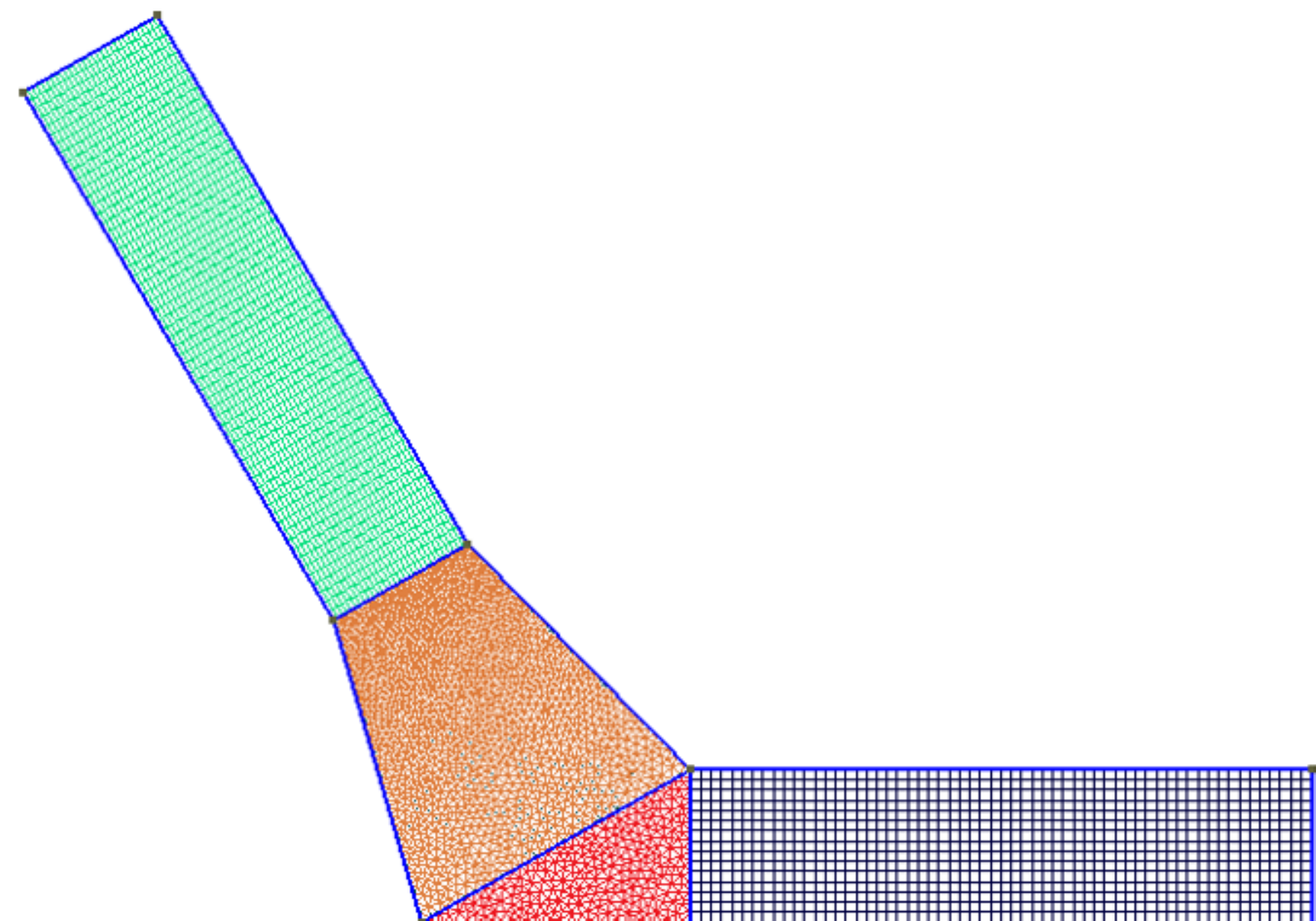}
	\end{subfigure}%
	\caption{(Left) Sketch of the geometry and mesh of Y-junction test case including boundary conditions. (Right) close up of the mesh in different zones.}
	\label{fig:Y_setup}
\end{figure}

Junctions are often used for the combination or separation of fluid flows and can be found in all types of engineering applications from gas transport in pipes till micro flow reactors. As a second test case a Y-junction with one outlet channel and two inlet channels is modeled. The angle between each inlet and the horizontal axis is 60 degrees, as shown in Figure~\ref{fig:Y_setup} on the left~\cite{Stephenson}. The length of the channels is 2 m. 

The 2D geometry is split in 6 zones as depicted in Figure~\ref{fig:Y_setup} on the left. On the three rectangular zones a mesh with quadrilateral cells is constructed. The remaining three zones are meshed with hexagonal cells. The different meshes are depicted in Figure~\ref{fig:Y_setup} on the right. The total number of cells is 13046. 

\begin{figure}[h!]
	\centering
	\begin{subfigure}{.42\textwidth}
		\centering
		\includegraphics[width=1.0\linewidth]{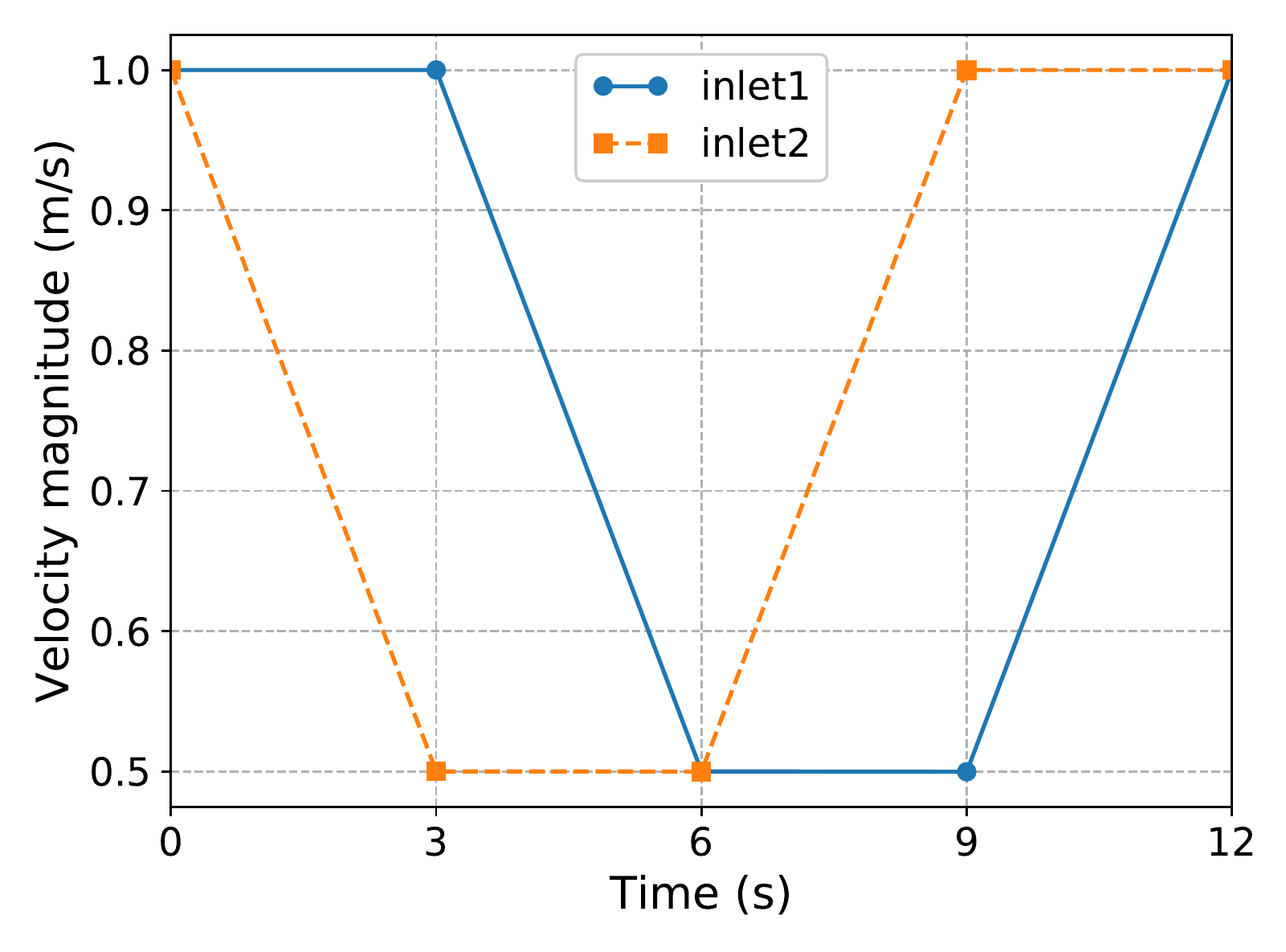}
	\end{subfigure}%
	\begin{subfigure}{.58\textwidth}
		\centering
		\includegraphics[width=1.0\linewidth]{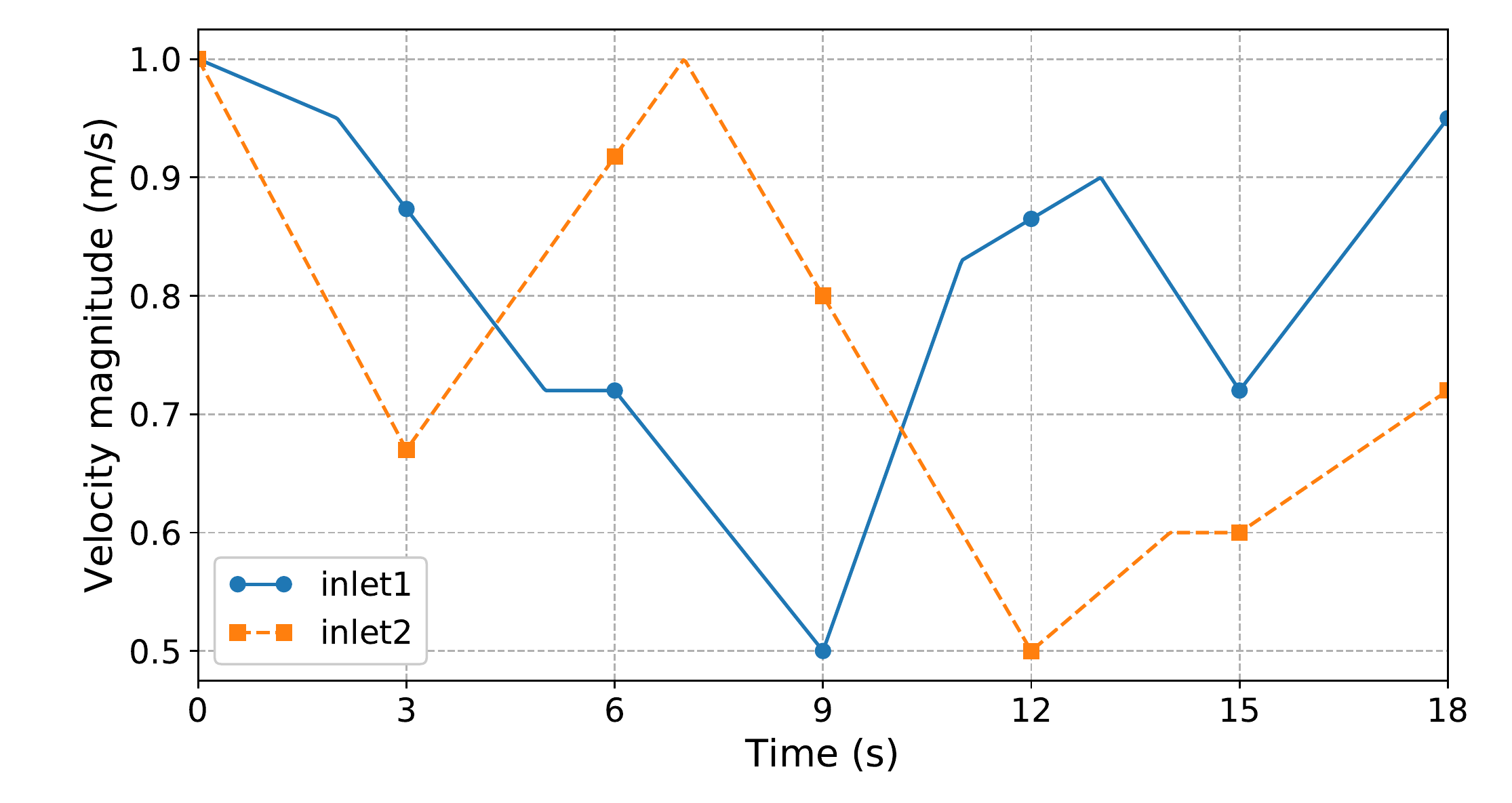}
	\end{subfigure}%
	\caption{Boundary conditions for the Y-junction test case. (Left) inlet velocity BCs for the FOM. (Right) inlet velocity BCs for ROM. }
	\label{fig:Y_BCs_eps}
\end{figure}

The boundary is subdivided into four different parts $\Gamma$ = $\Gamma_{i1}$ $\cup$ $\Gamma_{i2}$ $\cup$ $\Gamma_o$ $\cup$ $\Gamma_w$. The two inlets, $\Gamma_{i1}$ and $\Gamma_{i2}$, have a width of 0.5 m, while the outlet, $\Gamma_o$, has a width of 1 m. The kinematic viscosity is equal to $\nu$ = 1 $\cdot$ \num{e-2} m$^2$/s meaning that the Reynolds number at the inlet is 50 and the flow is considered laminar. The uniform inlet velocities are time-dependent and the velocity magnitude of the flow at the inlets is set according to figure~\ref{fig:Y_BCs_eps}.

A homogeneous Neumann boundary condition is applied for pressure at the inlet and wall boundaries. At the outlet, $\Gamma_o$, $p$ = 0 m$^2$/s$^2$ together with a homogeneous Neumann BC for velocity. A no slip BC is applied at the walls, $\Gamma_w$. 

As initial conditions the steady state solution, obtained with the $simpleFoam$ solver, for a velocity magnitude of 1 m/s at both inlets is chosen. The other boundary conditions are the same as for the unsteady simulation described above.

As done previously for the lid driven cavity case, the unsteady governing equations are iteratively solved by the FV method with the $pimpleFoam$ solver of OpenFOAM 6~\cite{Jasak}. For the full order simulations, the discretization in space is performed with a central differencing scheme for the diffusive term and a combination of a second order central-differencing and upwind schemes for the convective term. The temporal discretization is treated using a second order backward differencing scheme (BDF2). A constant time step of $\Delta t$ = 5 $\cdot$\num{e-4} s has been applied and the total full order simulation time is 12 s for which snapshots of the velocity and pressure fields are collected every 0.03 s, resulting in a total of 401 snapshots (including 1 for the initial condition). The inlet velocity BCs are time-dependent and the velocity magnitude of, alternately, inlet 1 or 2 is increased or decreased linearly between 1 m/s to 0.5 m/s as shown in Figure~\ref{fig:Y_BCs_eps} on the left.

\newpage
In that way, the ROM is trained for all possible combinations of inlet velocities within the specified range. The inlet boundary conditions of the ROM are then controlled according to Figure~\ref{fig:Y_BCs_eps} on the right, where the inlet velocity magnitude is increased or decreased linearly over time between the maximum of 1 m/s and minimum of 0.5 m/s. The magnitude of the inlet velocities of the ROM decreases and increases faster or slower over time compared to the training run. Also the ROM is tested for a longer time period, 18 s, compared to the full order simulation of 12 s. In that way, the ROM performance can be tested on the long term. 

The temporal discretization of the ROM is performed with a first order Newton's method. 

\begin{figure}[h!]
	\centering
	\begin{subfigure}{.35\textwidth}
		\centering
		\includegraphics[width=1.0\linewidth]{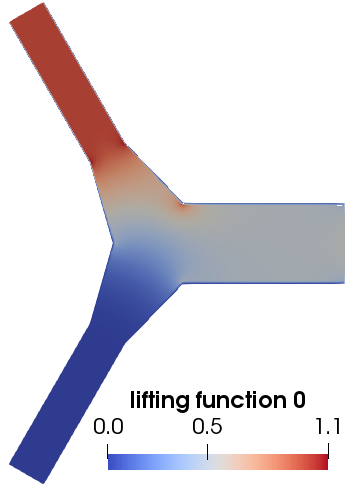}
	\end{subfigure}%
	\begin{subfigure}{.35\textwidth}
		\centering
		\includegraphics[width=1.0\linewidth]{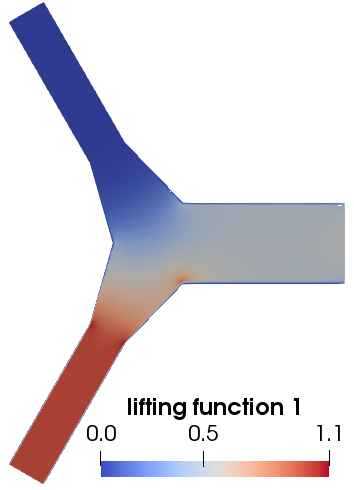}
	\end{subfigure}%
	\caption{The lifting functions for velocity for the Y-junction.}
	\label{fig:Y_control}
\end{figure}

Both the iterative penalty method and lifting function method are tested. The lifting functions are determined by solving for a potential flow field problem given by
\begin{align} \label{eq:potential_flow}
\begin{cases}
\nabla \cdot \boldsymbol{u} = 0   &\mbox{in  } \Omega,  \\
\nabla^2 p = 0 &\mbox{in  } \Omega,\\
\left(p(\boldsymbol{x})\boldsymbol{I}\right)\boldsymbol{n} = 0 &\mbox{on  } \Gamma_o, \\
\left(\nabla p(\boldsymbol{x},t) \right)\boldsymbol{n}= 0 &\mbox{on  }\Gamma \not\ni \Gamma_o, \\
\left(\nabla \boldsymbol{u}(\boldsymbol{x}) \right)\boldsymbol{n} = 0 &\mbox{on  } \Gamma_o,\\
\left(\nabla \boldsymbol{u}(\boldsymbol{x}) \right)\boldsymbol{n}  = 0 &\mbox{on  } \Gamma_{w}, \\
\boldsymbol{u}(\boldsymbol{x}) = \boldsymbol{g1}(\boldsymbol{x}) &\mbox{on  } \Gamma_{i1}, \\
\boldsymbol{u}(\boldsymbol{x}) = \boldsymbol{g2}(\boldsymbol{x}) &\mbox{on  } \Gamma_{i2}, 
\end{cases}
\end{align}
with the magnitude of the inlet velocity at inlet 1, $\Gamma_{i1}$, set to 1 m/s while inlet 2, $\Gamma_{i2}$, is kept at 0 m/s as shown in Figure ~\ref{fig:Y_control} for the first lifting function. To obtain the second lifting function $\|\boldsymbol{u}\|$ = 0 at $\Gamma_{i1}$ and 1 m/s at $\Gamma_{i2}$. Both lifting functions are shown in Figure~\ref{fig:Y_control}. 

The test case of a Y-junction is more complicated than the lid driven cavity case as not only one, but two boundaries need to be controlled, which are also time dependent. Furthermore, as the channel inlets are placed under an angle, one needs to take into account that the inlet velocity can be decomposed in an x- and a y-direction. Therefore, the vectorial lifting functions are split into their components before normalization. Also in the case of the penalty method, four penalty factors are determined; one for each inlet and each direction. This will be further discussed in Section~\ref{sec:discussion}.

\section{RESULTS AND ANALYSIS}\label{sec:results}
\subsection{Lid driven cavity flow problem}
First the full order simulation for the lid driven cavity test case is performed and 1001 velocity and pressure snapshots are collected, including the initial conditions, which are then used to create the POD basis functions. Stabile and Rozza~\cite{Stabile2017CAF} concluded in their research that 10 velocity and pressure modes are enough to retain 99.99$\%$ of the energy contained in the snapshots. Therefore, the same number of modes for the reduced basis creation are used in this work. 

Reduced order models are constructed with both the lifting function and penalty method and compared with a ROM without boundary enforcement. With the use of the iterative procedure a penalty factor of 0.058 is determined within 2 iterations by evaluating only the first five time steps with a maximum error tolerance, $\epsilon$, of~\num{e-5} for the value of the boundary condition of the ROM and starting from an initial guess of~\num{e-6}. For a similar study of the lid driven cavity benchmark problem, Lorenzi et al.~\cite{lorenzi2016pod} had found a factor between~\num{e-5} and~\num{e2} using numerical experimentation. The value found here using the iterative method is thus within the same range. A higher value for the penalty factor can be used, but it is then more likely that the ROM becomes ill-conditioned.

\begin{figure}[h!]
	\centering
	\begin{subfigure}{.49\textwidth}
		\centering
		\includegraphics[width=1.\linewidth]{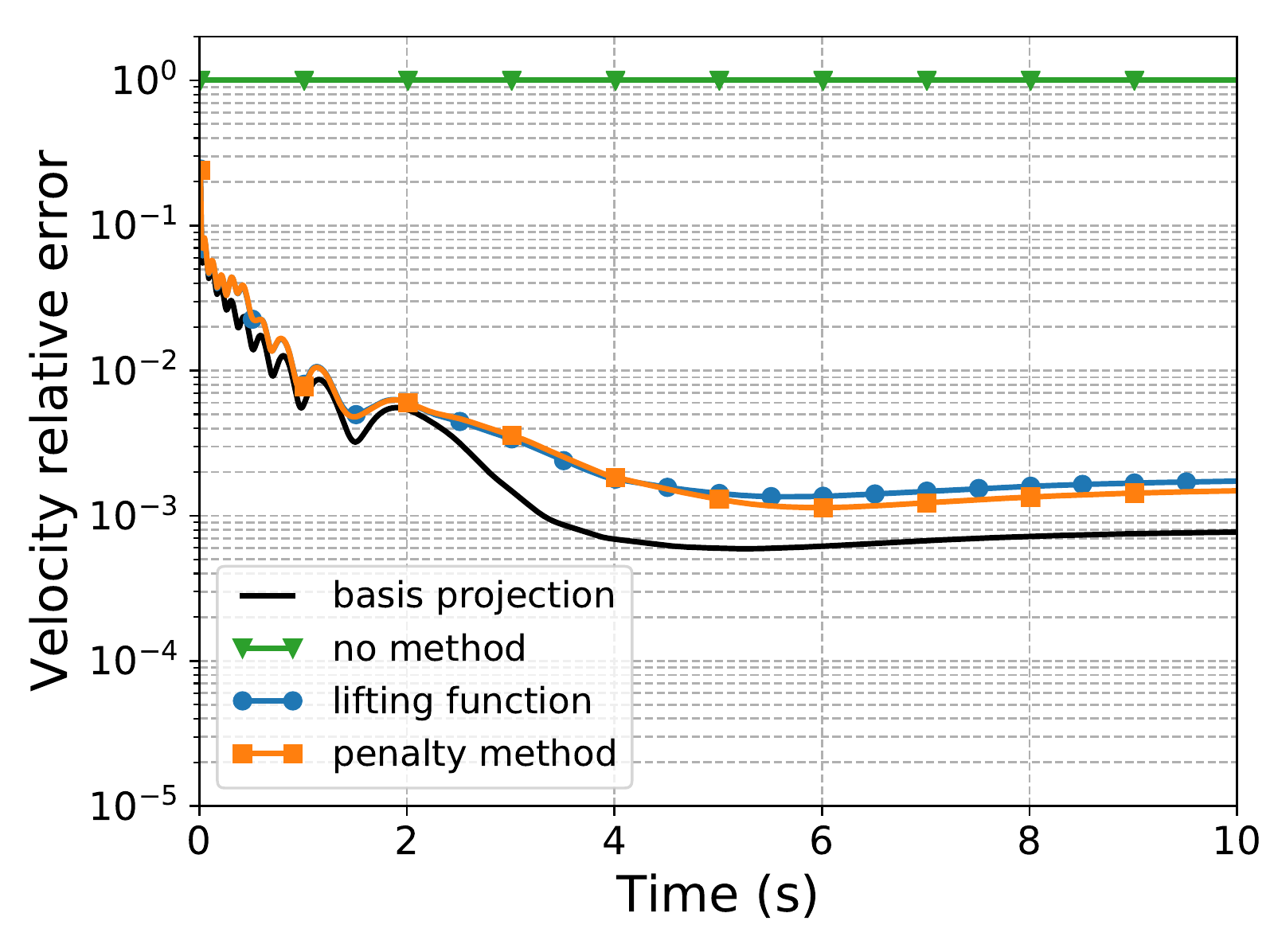}
	\end{subfigure}%
	\begin{subfigure}{.49\textwidth}
		\centering
		\includegraphics[width=1.\linewidth]{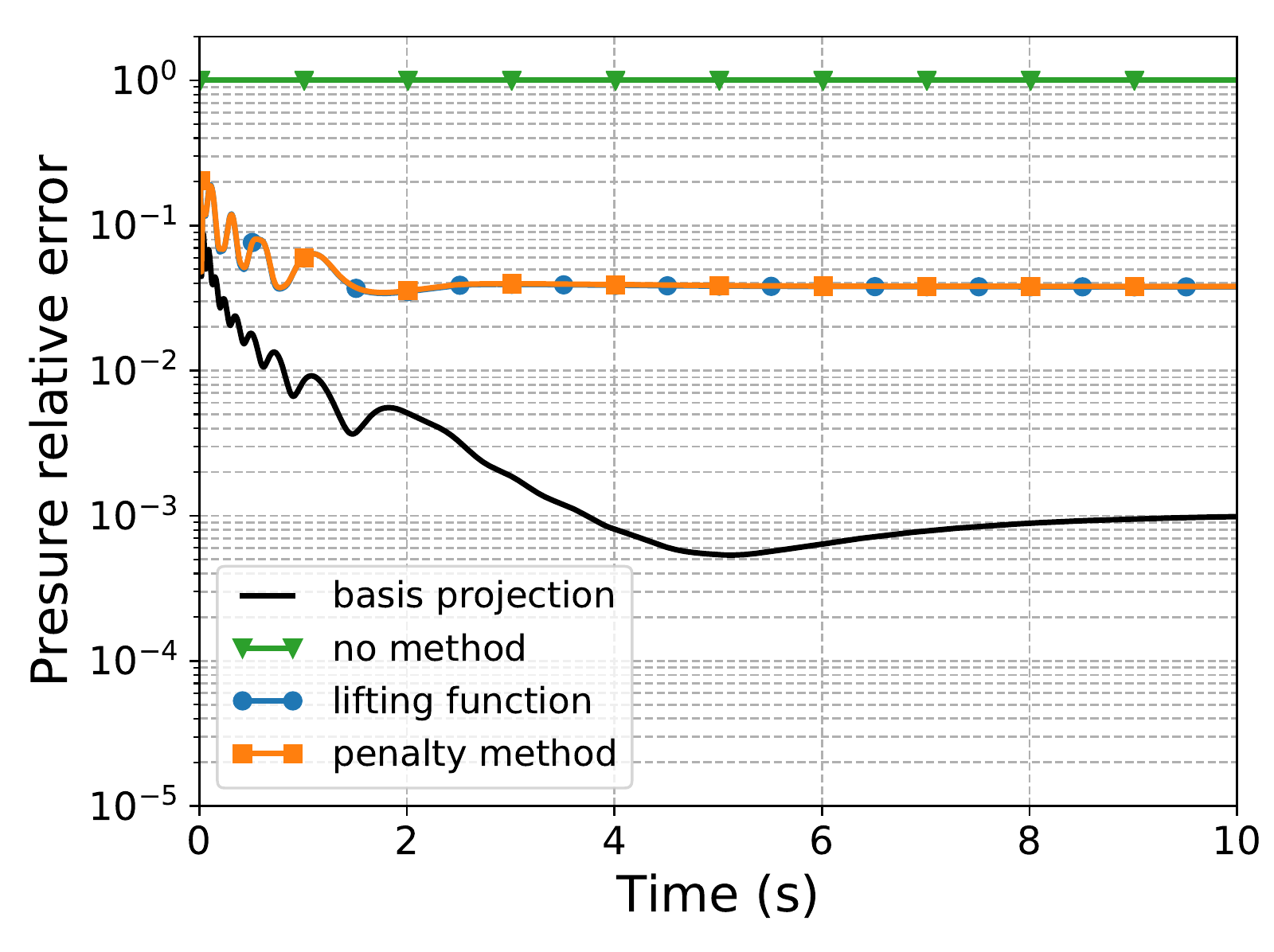}
	\end{subfigure}
	\caption{Relative $L^2$-error of velocity (left) and pressure (right) between the FOM and ROM with lifting function and with penalty method.}
	\label{fig:LID_L2_error}
\end{figure}

The obtained ROMs are tested for the same initial and boundary conditions as the high fidelity simulation. The evolution in time of the relative $L^2$-error between the reconstructed fields and the full order solutions is plotted in Figure~\ref{fig:LID_L2_error} together with the basis projection.

In case no boundary enforcement method is used the flow field remains zero throughout the simulation and therefore the relative error is 1. When either the lifting function or penalty method is used the relative $L^2$-error for both the velocity and pressure fields are about the order of \num{e-1} due to the relatively low number of snapshots acquired during the initial part of the transient. The snapshots are equally distributed in time, while this time span exhibits the most non-linear behavior. Therefore one should concentrate the snapshots in this time span to enhance the performance of the ROM~\cite{Stabile2017CAF}. After about 2 seconds of simulation time, for both boundary control methods, the velocity relative error drops till about the order of \num{e-3}. At the final time of the simulation the penalty method is performing slightly better than the lifting function method, but the order is the same.

\begin{figure}
	\centering
	\begin{subfigure}{.19\linewidth}
		\centering
		\includegraphics[width=1.\linewidth]{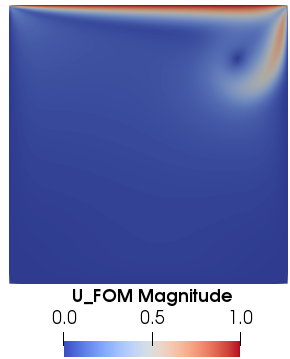}
	\end{subfigure}%
	\begin{subfigure}{.19\textwidth}
		\centering
		\includegraphics[width=1.\linewidth]{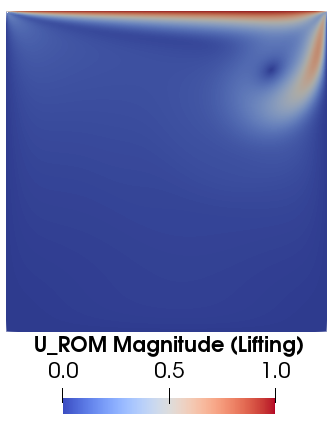}
	\end{subfigure}
	\begin{subfigure}{.19\textwidth}
		\centering
		\includegraphics[width=1.\linewidth]{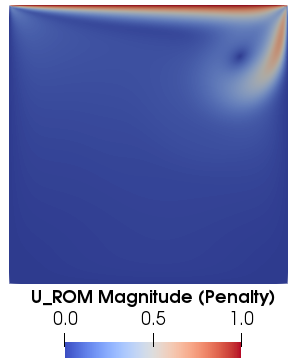}
	\end{subfigure}
	\begin{subfigure}{.19\textwidth}
		\centering
		\includegraphics[width=1.\linewidth]{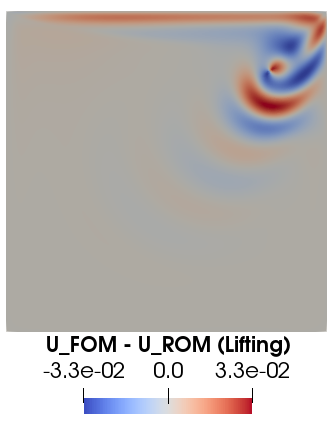}
	\end{subfigure}
	\begin{subfigure}{.19\textwidth}
		\centering
		\includegraphics[width=1.\linewidth]{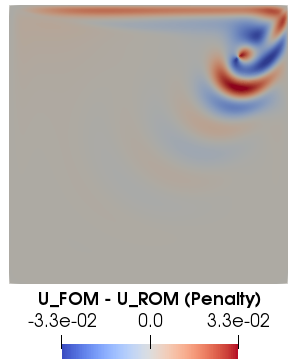}
	\end{subfigure}
	
	\begin{subfigure}{.19\linewidth}
		\centering
		\includegraphics[width=1.\linewidth]{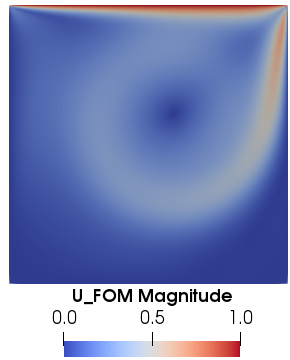}
	\end{subfigure}%
	\begin{subfigure}{.19\textwidth}
		\centering
		\includegraphics[width=1.\linewidth]{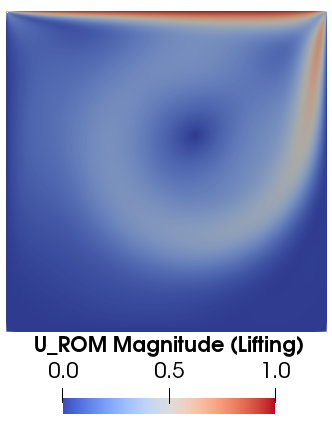}
	\end{subfigure}
	\begin{subfigure}{.19\textwidth}
		\centering
		\includegraphics[width=1.\linewidth]{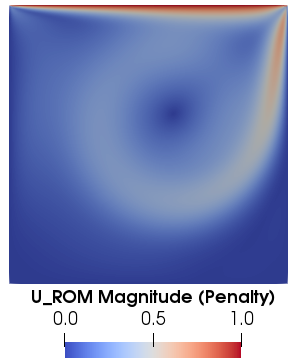}
	\end{subfigure}
	\begin{subfigure}{.19\textwidth}
		\centering
		\includegraphics[width=1.\linewidth]{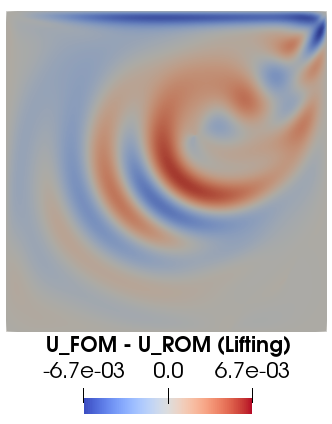}
	\end{subfigure}
	\begin{subfigure}{.19\textwidth}
		\centering
		\includegraphics[width=1.\linewidth]{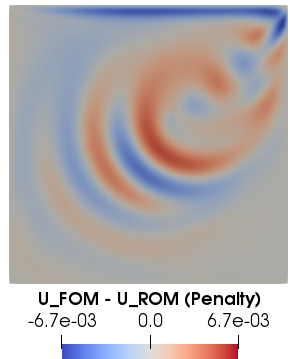}
	\end{subfigure}
	
	\begin{subfigure}{.19\linewidth}
		\centering
		\includegraphics[width=1.\linewidth]{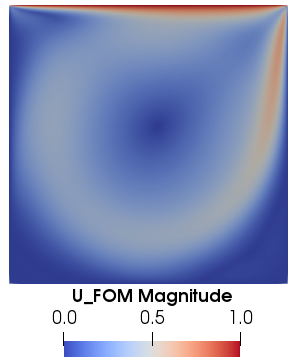}
	\end{subfigure}%
	\begin{subfigure}{.19\textwidth}
		\centering
		\includegraphics[width=1.\linewidth]{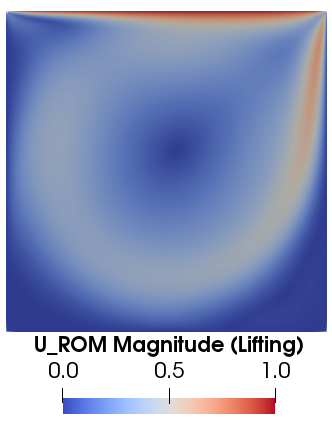}
	\end{subfigure}
	\begin{subfigure}{.19\textwidth}
		\centering
		\includegraphics[width=1.\linewidth]{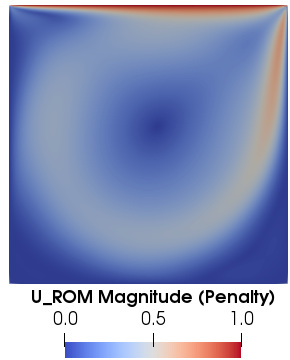}
	\end{subfigure}
	\begin{subfigure}{.19\textwidth}
		\centering
		\includegraphics[width=1.\linewidth]{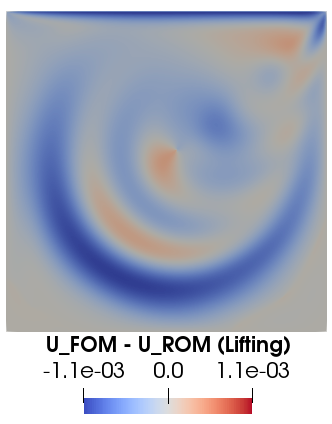}
	\end{subfigure}
	\begin{subfigure}{.19\textwidth}
		\centering
		\includegraphics[width=1.\linewidth]{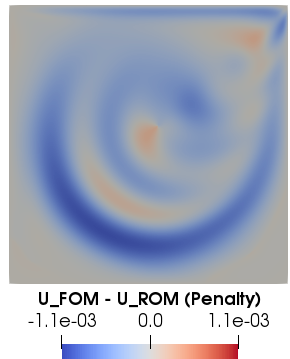}
	\end{subfigure}
	
	\begin{subfigure}{.19\linewidth}
		\centering
		\includegraphics[width=1.\linewidth]{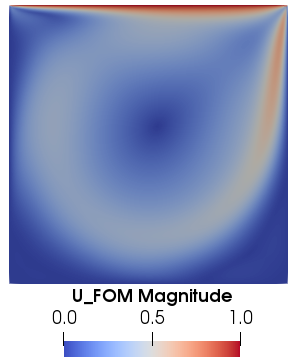}
	\end{subfigure}%
	\begin{subfigure}{.19\textwidth}
		\centering
		\includegraphics[width=1.\linewidth]{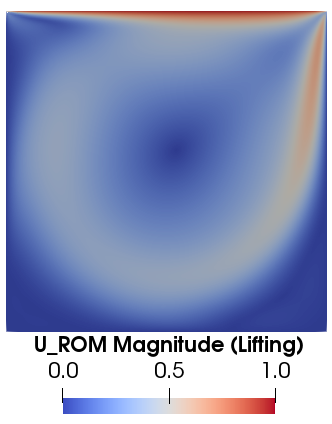}
	\end{subfigure}
	\begin{subfigure}{.19\textwidth}
		\centering
		\includegraphics[width=1.\linewidth]{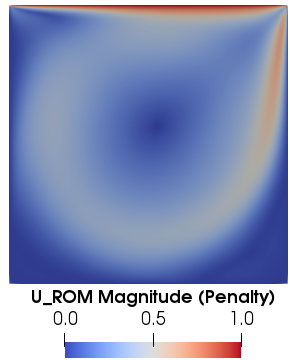}
	\end{subfigure}
	\begin{subfigure}{.19\textwidth}
		\centering
		\includegraphics[width=1.\linewidth]{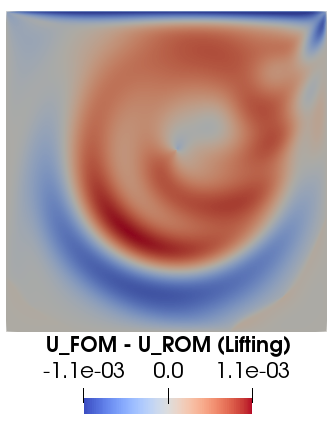}
	\end{subfigure}
	\begin{subfigure}{.19\textwidth}
		\centering
		\includegraphics[width=1.\linewidth]{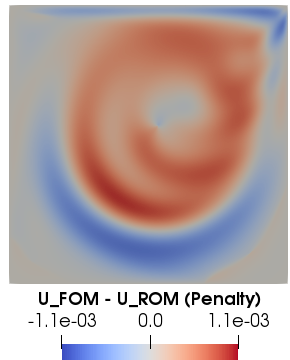}
	\end{subfigure}
	\caption{Comparison of the full order velocity magnitude fields (1st column), the ROM fields obtained with the lifting function method (2nd column) and penalty method (4th column) and the difference between the FOM and ROM fields obtained with the lifting function method (3rd column) and penalty method (5th column) at $t$ = 0.2, 1, 5 and 10 s (from top to bottom) for the lid driven cavity problem.}
	\label{fig:LID_U}
\end{figure}

\begin{figure}
	\centering
	\begin{subfigure}{.19\linewidth}
		\centering
		\includegraphics[width=1.\linewidth]{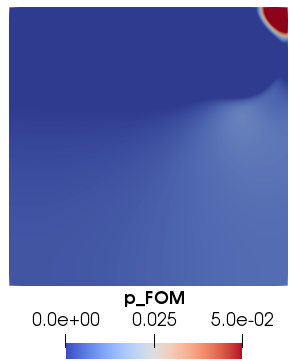}
	\end{subfigure}%
	\begin{subfigure}{.19\textwidth}
		\centering
		\includegraphics[width=1.\linewidth]{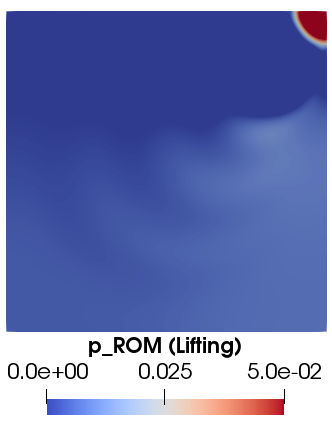}
	\end{subfigure}
	\begin{subfigure}{.19\textwidth}
		\centering
		\includegraphics[width=1.\linewidth]{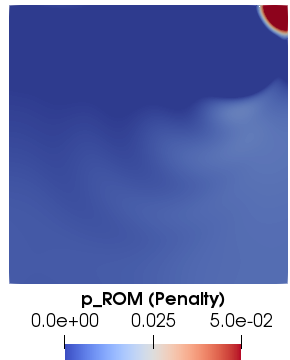}
	\end{subfigure}
	\begin{subfigure}{.19\textwidth}
		\centering
		\includegraphics[width=1.\linewidth]{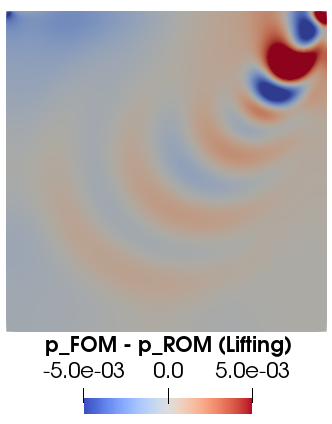}
	\end{subfigure}
	\begin{subfigure}{.19\textwidth}
		\centering
		\includegraphics[width=1.\linewidth]{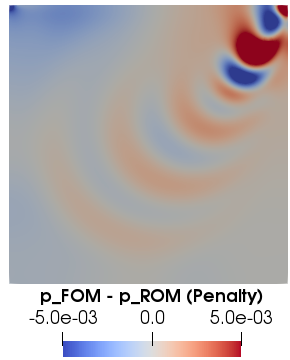}
	\end{subfigure}
	
	\begin{subfigure}{.19\linewidth}
		\centering
		\includegraphics[width=1.\linewidth]{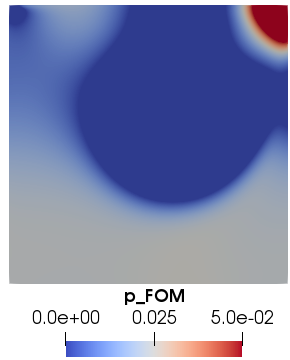}
	\end{subfigure}%
	\begin{subfigure}{.19\textwidth}
		\centering
		\includegraphics[width=1.\linewidth]{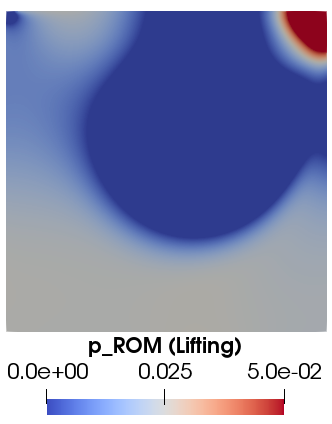}
	\end{subfigure}
	\begin{subfigure}{.19\textwidth}
		\centering
		\includegraphics[width=1.\linewidth]{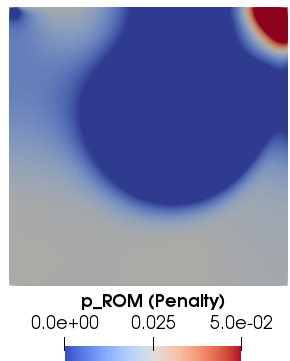}
	\end{subfigure}
	\begin{subfigure}{.19\textwidth}
		\centering
		\includegraphics[width=1.\linewidth]{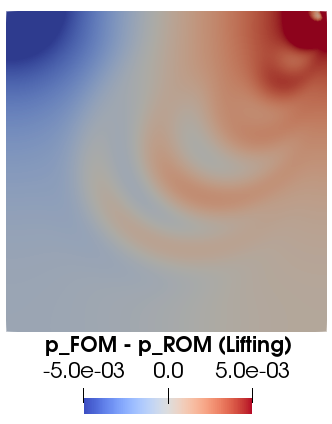}
	\end{subfigure}
	\begin{subfigure}{.19\textwidth}
		\centering
		\includegraphics[width=1.\linewidth]{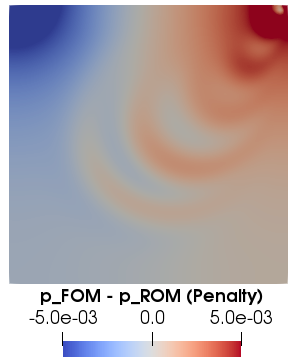}
	\end{subfigure}
	
	\begin{subfigure}{.19\linewidth}
		\centering
		\includegraphics[width=1.\linewidth]{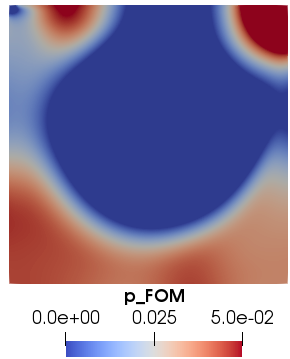}
	\end{subfigure}%
	\begin{subfigure}{.19\textwidth}
		\centering
		\includegraphics[width=1.\linewidth]{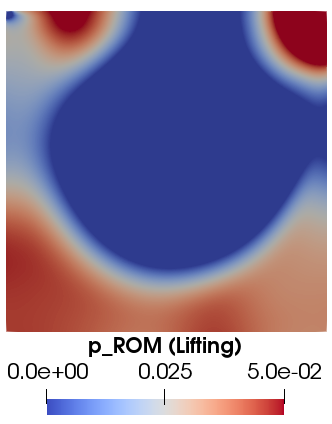}
	\end{subfigure}
	\begin{subfigure}{.19\textwidth}
		\centering
		\includegraphics[width=1.\linewidth]{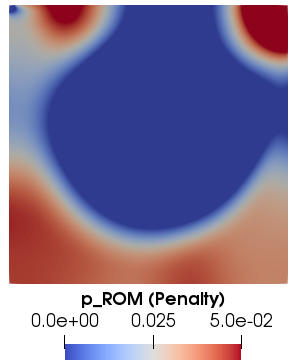}
	\end{subfigure}
	\begin{subfigure}{.19\textwidth}
		\centering
		\includegraphics[width=1.\linewidth]{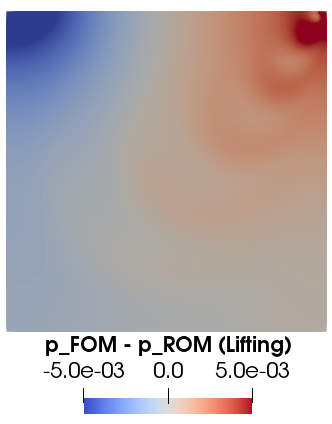}
	\end{subfigure}
	\begin{subfigure}{.19\textwidth}
		\centering
		\includegraphics[width=1.\linewidth]{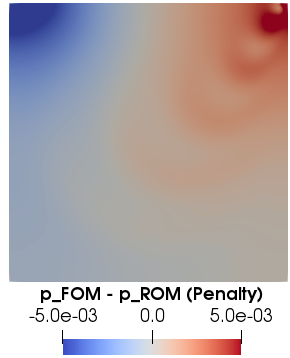}
	\end{subfigure}
	
	\begin{subfigure}{.19\linewidth}
		\centering
		\includegraphics[width=1.\linewidth]{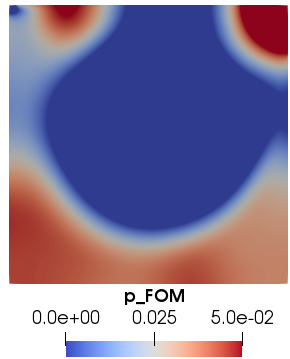}
	\end{subfigure}%
	\begin{subfigure}{.19\textwidth}
		\centering
		\includegraphics[width=1.\linewidth]{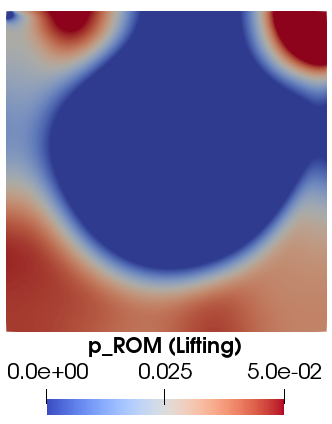}
	\end{subfigure}
	\begin{subfigure}{.19\textwidth}
		\centering
		\includegraphics[width=1.\linewidth]{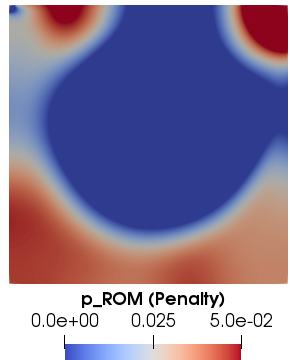}
	\end{subfigure}
	\begin{subfigure}{.19\textwidth}
		\centering
		\includegraphics[width=1.\linewidth]{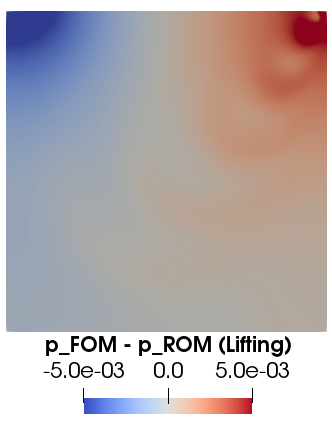}
	\end{subfigure}
	\begin{subfigure}{.19\textwidth}
		\centering
		\includegraphics[width=1.\linewidth]{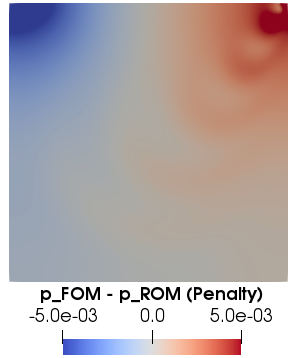}
	\end{subfigure}
	\caption{Comparison of the full order pressure fields (1st column), the ROM fields obtained with the lifting function method (2nd column) and penalty method (4th column) and the difference between the FOM and ROM fields obtained with the lifting function method (3rd column) and penalty method (5th column) at $t$ = 0.2, 1, 5 and 10 s (from top to bottom) for the lid driven cavity problem. }
	\label{fig:LID_p}
\end{figure}

Contrary to velocity, the relative error for pressure stays about 4$\cdot$\num{e-1} after 2 s of simulation time, while the projection error drops till about \num{e-3}. This has been previously acclaimed by Stabile et al. in~\cite{Stabile2017CAIM}. The PPE stabilization method is less accurate concerning pressure compared to the supremizer enrichment method. This has also been found by Kean and Schneier~\cite{kean2020error} in the finite element-based ROM setting.
Furthermore, the absolute error between the FOM and the ROMs is shown in Figure~\ref{fig:LID_U} and~\ref{fig:LID_p} for velocity magnitude and pressure, respectively. 

It is observed that both methods lead, for velocity, to an absolute error between the FOM and the ROM of the order \num{e-2} at the beginning of the simulation and about \num{e-3} once the flow has reached its steady state solution. Furthermore, the velocity error slightly increases between 5 and 10 s of simulation time. This can also be observed in the $L^2$-error analysis over time in Figure~\ref{fig:LID_L2_error}. For pressure, the error is largest near the top corners of the cavity and are of the order \num{e-3}. Note that the scale does not show the whole range of absolute errors. This is done to better visualize the error. The maximum error for pressure is about 5$\cdot$\num{e-2} m$^2$/s$^2$ at the top right corner. As the pressure relative to its reference point at (0,0) plotted in Figure~\ref{fig:LID_p} is always less than 1  m$^2$/s$^2$, the relative error plotted in Figure~\ref{fig:LID_L2_error} is greater than the absolute error plotted in Figure~\ref{fig:LID_p}. Furthermore, the error distribution, for both the velocity and pressure fields, is similar all over the domain, meaning the methods are performing the same, as previously confirmed by the $L^2$-error analysis over time in Figure~\ref{fig:LID_L2_error}.

\begin{figure}[t!]
	\centering
	\includegraphics[width=0.5\linewidth]{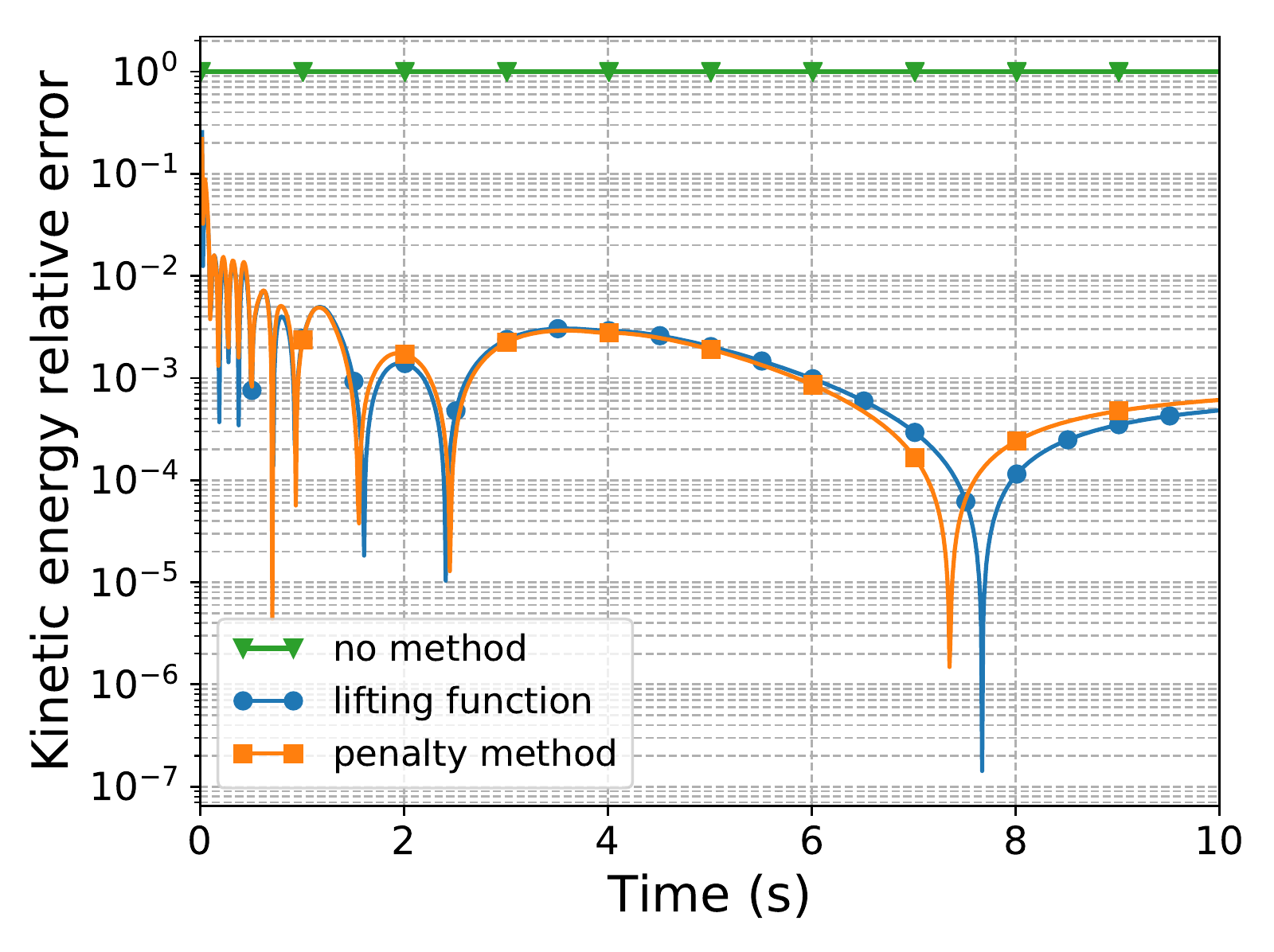}
	\caption{Kinetic energy relative $L^2$-error for the ROM with lifting function and with penalty method.}
	\label{fig:LID_KE}
\end{figure}

The relative error for the total kinetic energy is determined and plotted in Figure~\ref{fig:LID_KE}. The order is more or less the same for both boundary control methods. From time to time the penalty method is performing slightly better and the other way around and the relative velocity error is less than \num{e-2} for the vast part of the simulation.

\newpage
Finally, the computational times for performing the full order simulation (Eq.~\ref{eq:FOM_mat}), calculating the POD modes (Eqs.~\ref{eq:approx}-~\ref{eq:POD}), the reduced matrices (Eqs.~\ref{eq:ROM_matrices},~\ref{eq:C_matrix} and~\ref{eq:ROM_matrices2}) and performing the simulation at reduced level (Eq.~\ref{eq:ROM} (lifting function method) or Eq.~\ref{eq:pen_ROM_min} (penalty method) \& Eq.~\ref{eq:ROM2}) are all listed in Table~\ref{tab:times}. Calculating the POD modes, reduced matrices and the ROM solutions takes more time in the case of the lifting function method as the reduced basis space consists of an additional mode, namely the normalized lifting function for the boundary with the lid, compared to the penalty method. Determining the penalty factor with the iterative method takes only 0.11 s. The speedup ratio between the ROM and the FOM is about 270 times for the lifting method and 308 times for the penalty method.

\begin{table}[h!]
	\caption{Computational time (clock time) for the FOM simulation, POD, calculating reduced matrices offline (Matrices), determining penalty factor with iterative method (Penalty factor) and ROM simulation.}
	\centering
	\begin{tabular}{lllllll}
		\hline
		\multicolumn{1}{c}{Method}  &FOM  &POD  &Matrices  &Penalty factor &ROM  \\ \hline
		Lifting                       &37 min.     & 50 s  & 8.2 s & - &8.2 s    \\
		Penalty 	  & 37 min. & 45 s  & 6.8 s & 0.11 s &7.2 s  \\ \hline
	\end{tabular}
	\label{tab:times}
\end{table}


\clearpage
\subsection{Y-junction flow problem}
A full order simulation is performed for the Y-junction test case with varying inlet velocities (magnitude) according to figure~\ref{fig:Y_BCs_eps} on the left. In total 401 velocity and pressure snapshots are collected, which are then used to created the POD basis functions. To determine the number of basis functions necessary for the creation of the reduced subspace, the cumulative eigenvalues (based on the first 20 most energetic POD modes) are listed in Table~\ref{tab:Y_cumm_ev}. 

\begin{figure}[h!]
	\centering
	\begin{subfigure}{.195\textwidth}
		\centering
		\includegraphics[width=1.\linewidth]{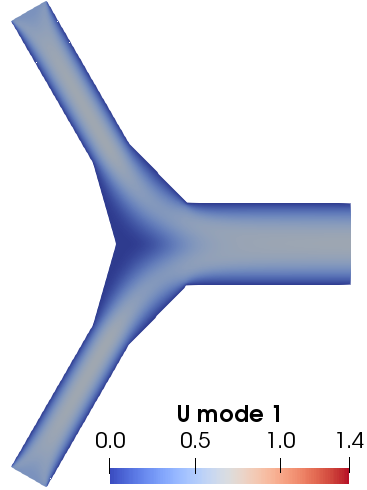}
	\end{subfigure}%
	\begin{subfigure}{.195\textwidth}
		\centering
		\includegraphics[width=1.\linewidth]{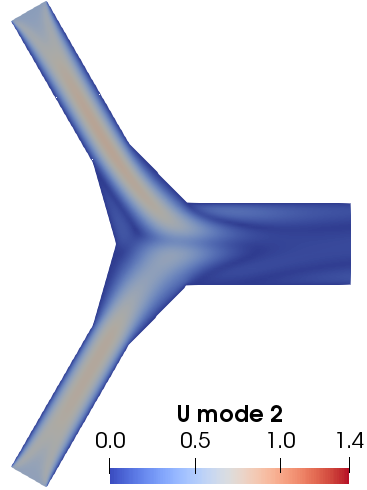}
	\end{subfigure}
	\begin{subfigure}{.195\textwidth}
		\centering
		\includegraphics[width=1.\linewidth]{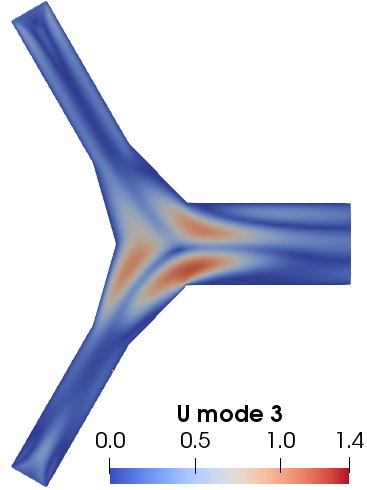}
	\end{subfigure}
	\begin{subfigure}{.195\textwidth}
		\centering
		\includegraphics[width=1.\linewidth]{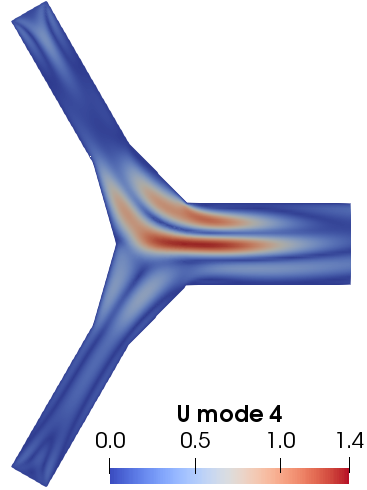}
	\end{subfigure}
	\begin{subfigure}{.195\textwidth}
		\centering
		\includegraphics[width=1.\linewidth]{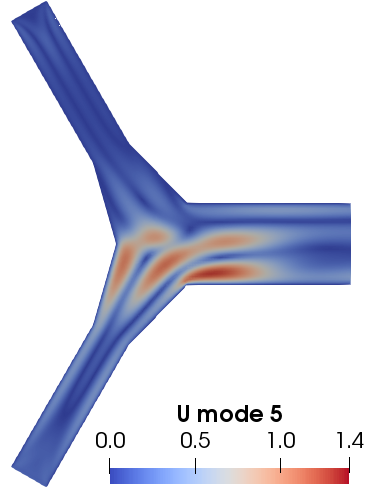}
	\end{subfigure}
	\begin{subfigure}{.195\textwidth}
		\centering
		\includegraphics[width=1.\linewidth]{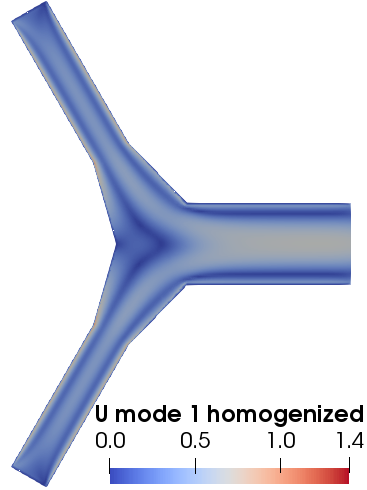}
	\end{subfigure}%
	\begin{subfigure}{.195\textwidth}
		\centering
		\includegraphics[width=1.\linewidth]{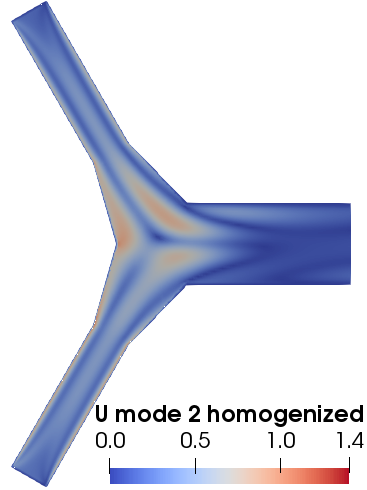}
	\end{subfigure}
	\begin{subfigure}{.195\textwidth}
		\centering
		\includegraphics[width=1.\linewidth]{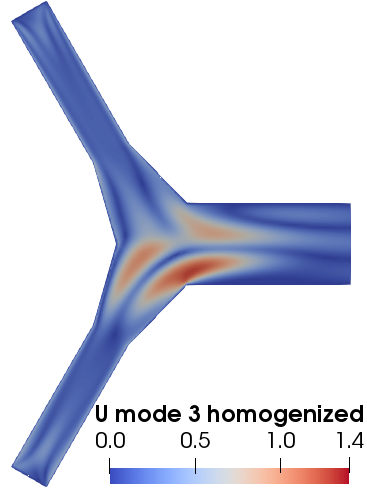}
	\end{subfigure}
	\begin{subfigure}{.195\textwidth}
		\centering
		\includegraphics[width=1.\linewidth]{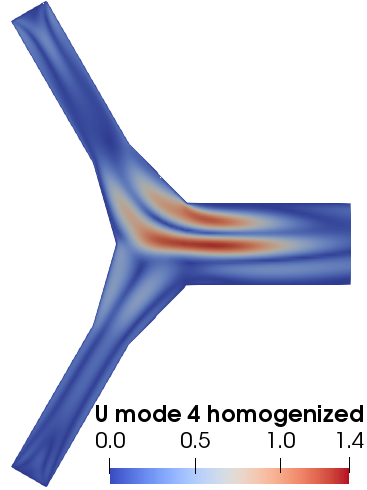}
	\end{subfigure}
	\begin{subfigure}{.195\textwidth}
		\centering
		\includegraphics[width=1.\linewidth]{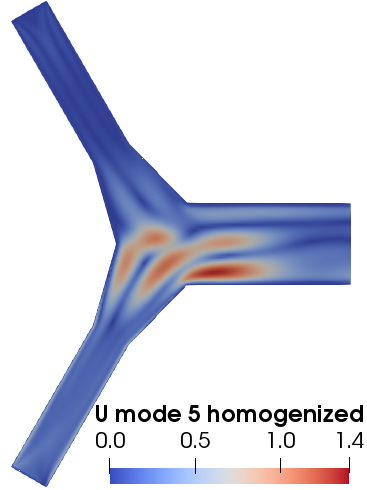}
	\end{subfigure}
	\begin{subfigure}{.195\textwidth}
		%\centering
		\includegraphics[width=1.\linewidth]{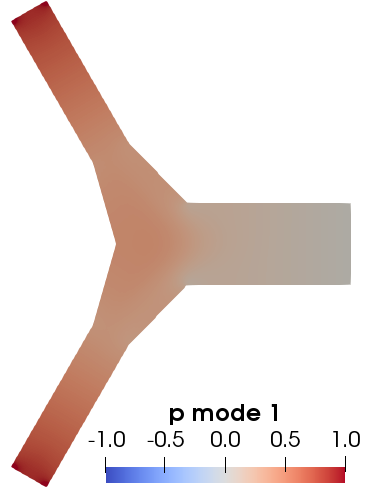}
	\end{subfigure}%
	\begin{subfigure}{.195\textwidth}
		%\centering
		\includegraphics[width=1.\linewidth]{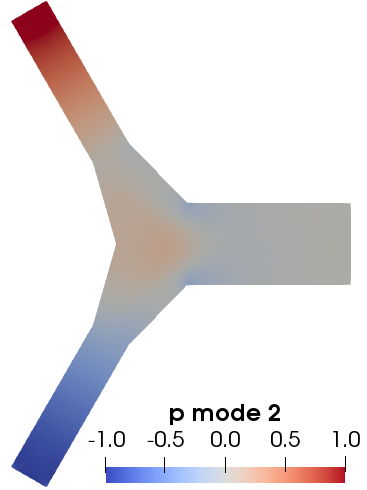}
	\end{subfigure}
	\begin{subfigure}{.195\textwidth}
		%\centering
		\includegraphics[width=1.\linewidth]{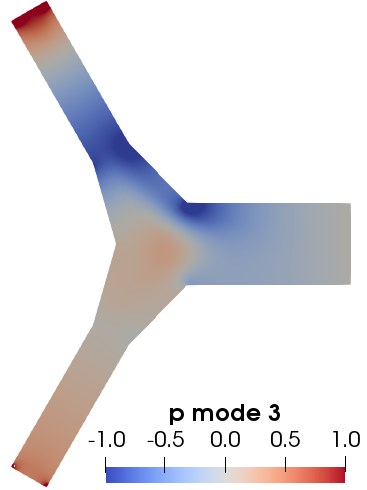}
	\end{subfigure}
	\begin{subfigure}{.195\textwidth}
		%\centering
		\includegraphics[width=1.\linewidth]{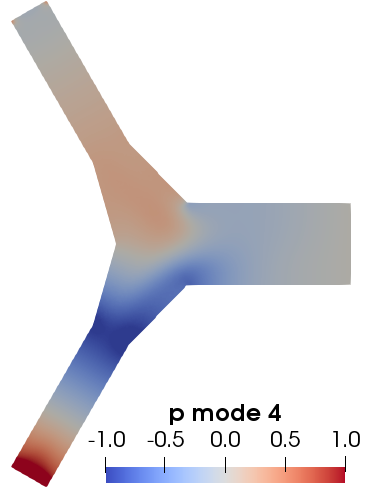}
	\end{subfigure}
	\begin{subfigure}{.195\textwidth}
		%	\centering
		\includegraphics[width=1.\linewidth]{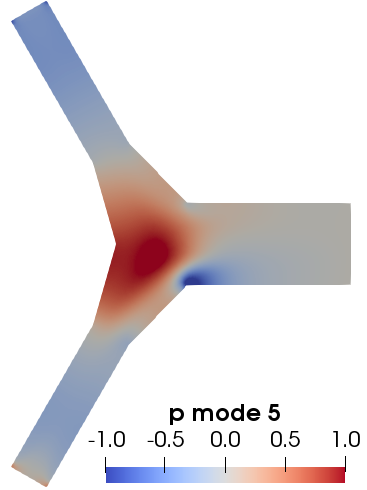}
	\end{subfigure}
	\caption{First 5 POD modes for (top) velocity, (middle) velocity with homogeneous BCs and (bottom) pressure for the Y-junction flow problem. }
	\label{fig:Y_modes}
\end{figure}

\begin{table}[h!]
	\caption{The cumulative eigenvalues for the Y-junction test case. The second and third columns report the cumulative eigenvalues (total of the first 20 modes) for the velocity and pressure fields, respectively.}
	\centering
	\begin{tabular}{lll}
		\hline
		\multicolumn{1}{l}{N modes} & \textbf{u} & $p$  \\ \hline
		1                           & 0.976478 & 0.967073  \\
		2                           & 0.998492 & 0.989840  \\
		3                           & 0.999724 & 0.998781  \\
		4                           & 0.999859 & 0.999741  \\
		5                           & 0.999924 & 0.999933  \\
		6							& 0.999967 & 0.999975  \\
		7 							& 0.999989 & 0.999995 \\
		10     				    	& 0.999999 & 0.999999 \\ \hline
	\end{tabular}
	\label{tab:Y_cumm_ev}
\end{table}
5 velocity and pressure modes are sufficient to retain 99.99$\%$ of the energy contained in the snapshots. These first five (homogenized) velocity and pressure modes are plotted in Figure~\ref{fig:Y_modes}. The first velocity magnitude mode has a symmetric pattern and is close to the time-averaged solution when it has non-homogeneous BCs and looks more like a fluctuation around the mean when it has homogeneous BCs. From the third mode and higher, the modes are more or less alike, whether the modes have homogeneous BCs or not.

\begin{figure}[h!]
	\centering
	\includegraphics[width=0.95\linewidth]{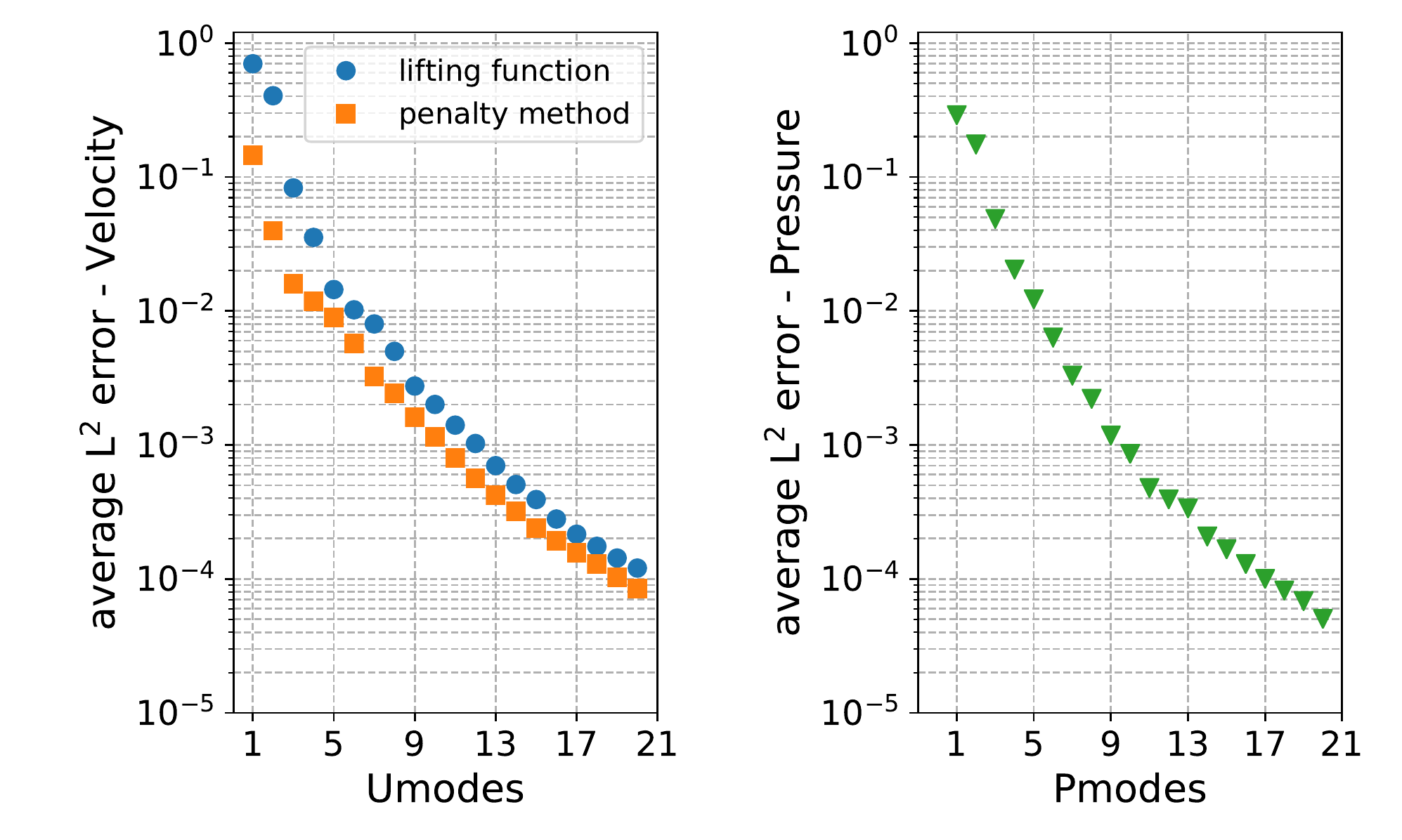}
	\caption{The time-averaged $L^2$-error per number of (left) velocity modes (Umodes) and (right) pressure modes (Pmodes) for the Y-junction test case.}
	\label{fig:Y_ev}
\end{figure}

In Figure~\ref{fig:Y_ev} for each number of modes the time-averaged relative $L^2$-error between the FOM and the basis projection is plotted, on the left for velocity and on the right for pressure. For velocity this is repeated with a set of homogenized snapshots. As there are two inlet boundary conditions, the first two modes are the normalized lifting functions and all sequential modes are then the homogeneous basis functions obtained with the POD method. Therefore the average $L^2$-error is still above the order \num{e-1} as these modes do not contain any information about the full order solution. The figure shows that 11 velocity basis functions and 10 pressure basis functions are required to have a truncation error less than \num{e-3}. Taking also into account previous observation, these number of modes are used for calculating the ROM matrices. \newpage

After applying the Galerkin projection with the obtained modes, the penalty factors are determined using the iterative procedure. Starting from an initial guess of \num{e-6} the penalty factors found are 5.9$\cdot$\num{e-8} and 88.3 for inlet 1 and 1.1$\cdot$\num{e-7} and 125 for inlet 2 in the x-direction and y-direction, respectively. The factors are determined within 41 iterations for an error tolerance of \num{e-5} and only the first five time steps are evaluated. However, it took only 15 iterations to have an error of 1.00009$\cdot$\num{e-5} with penalty factors 0.0327, 88.3, 0.048, 124.5. So one could relax the criteria for the error a bit for a faster convergence. \newpage

\begin{figure}[h!]
	\centering
	\begin{subfigure}{.49\textwidth}
		\centering
		\includegraphics[width=1.\linewidth]{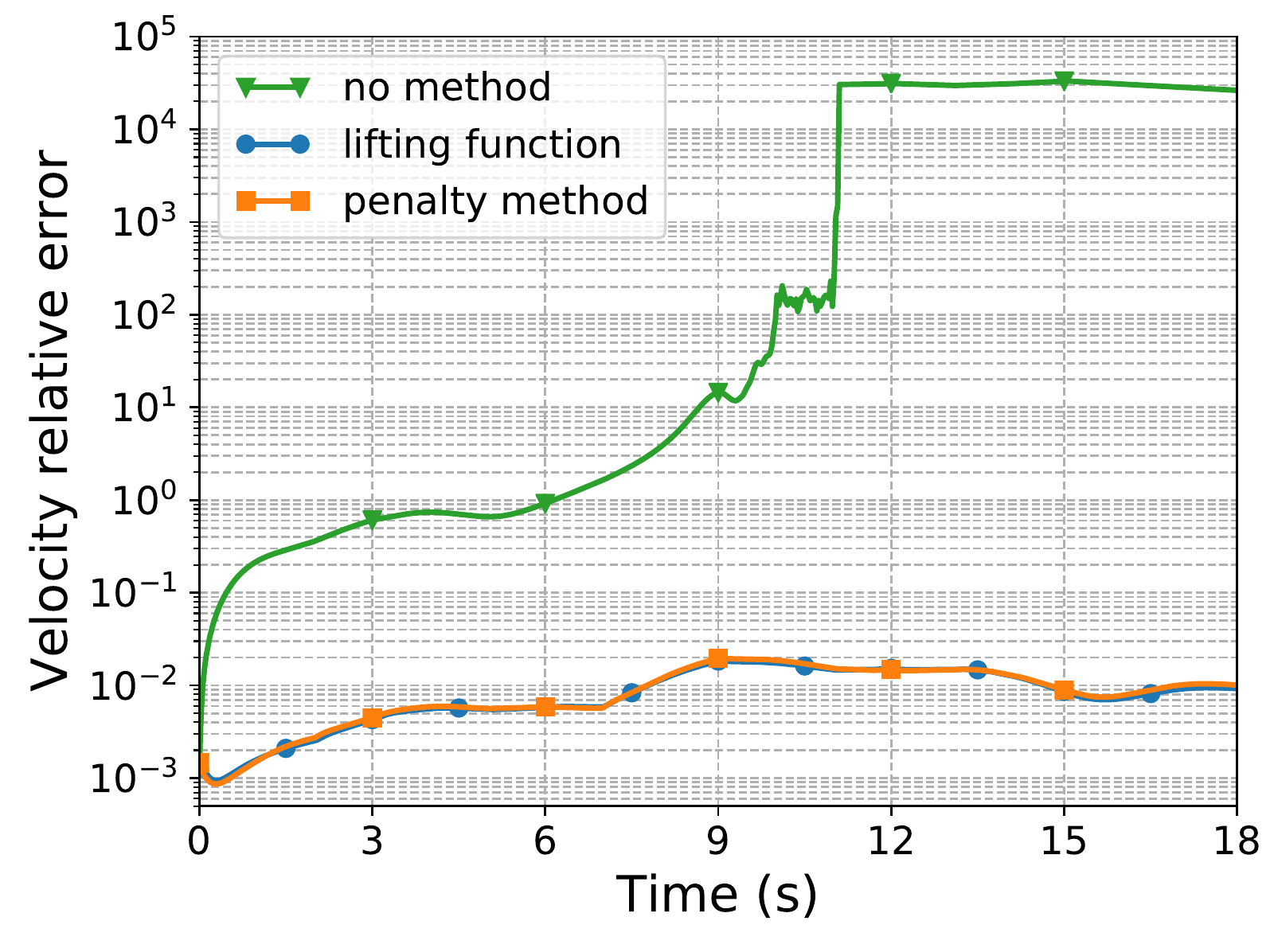}
	\end{subfigure}%
	\begin{subfigure}{.49\textwidth}
		\centering
		\includegraphics[width=1.\linewidth]{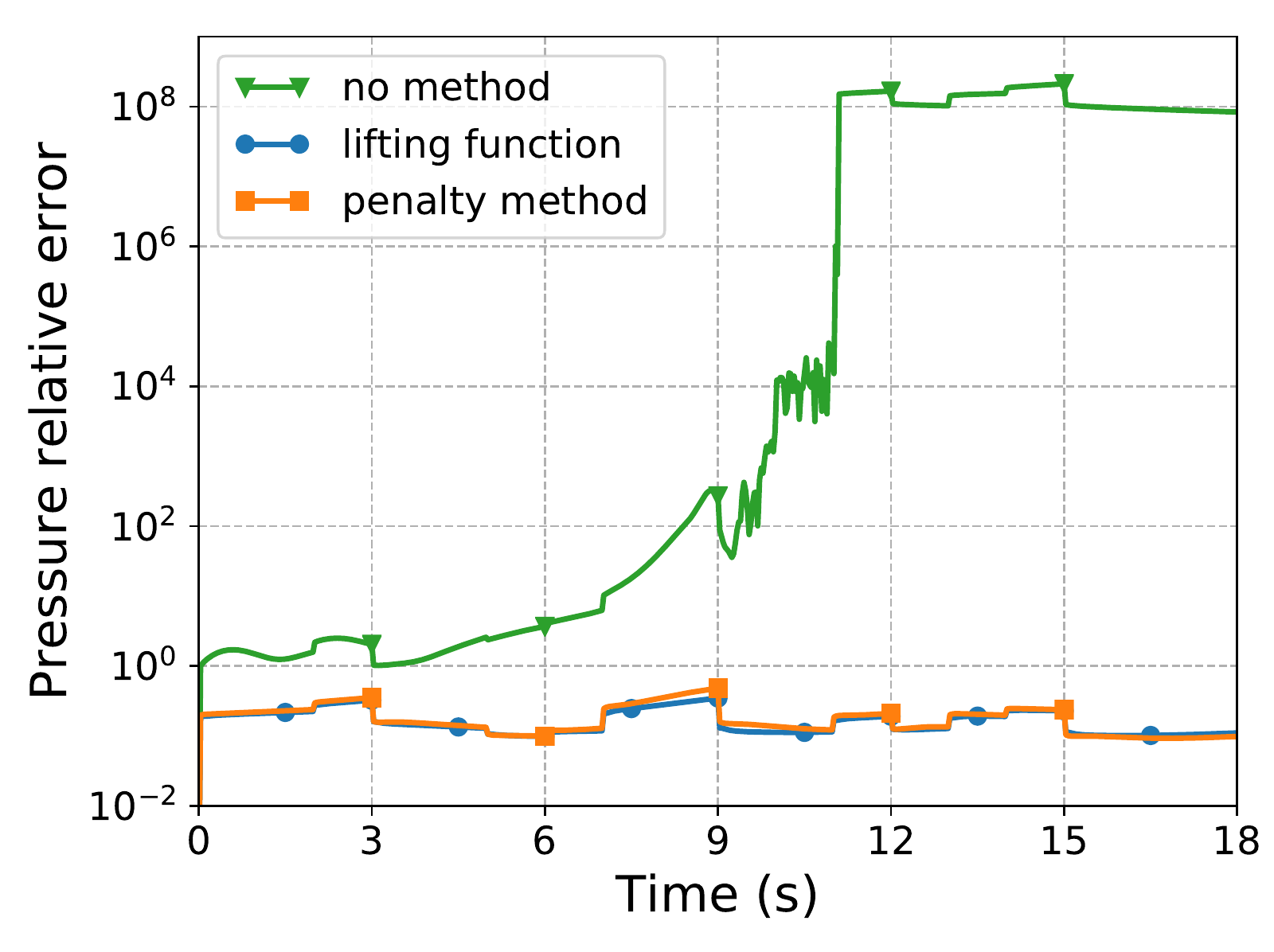}
	\end{subfigure}
	\caption{Relative $L^2$-error of velocity (left) and pressure (right) between the FOM and ROM with lifting function and with penalty method for the Y-junction flow problem.}
	\label{fig:Y_L2_error}
\end{figure}

Thereafter, three ROMs are obtained; one without boundary enforcement method, one with the lifting function method and one with the penalty method. These are then consecutively tested for the time-dependent boundary conditions of Figure~\ref{fig:Y_control}. The evolution in time of the $L^2$ relative error between the reconstructed fields is plotted in Figure~\ref{fig:Y_L2_error}.

In case no boundary enforcement method is used, the relative error for both velocity and pressure is of the order 1 and larger for the vast part of the simulation. 

\begin{figure}[h!]
	\centering
	\includegraphics[width=0.50\linewidth]{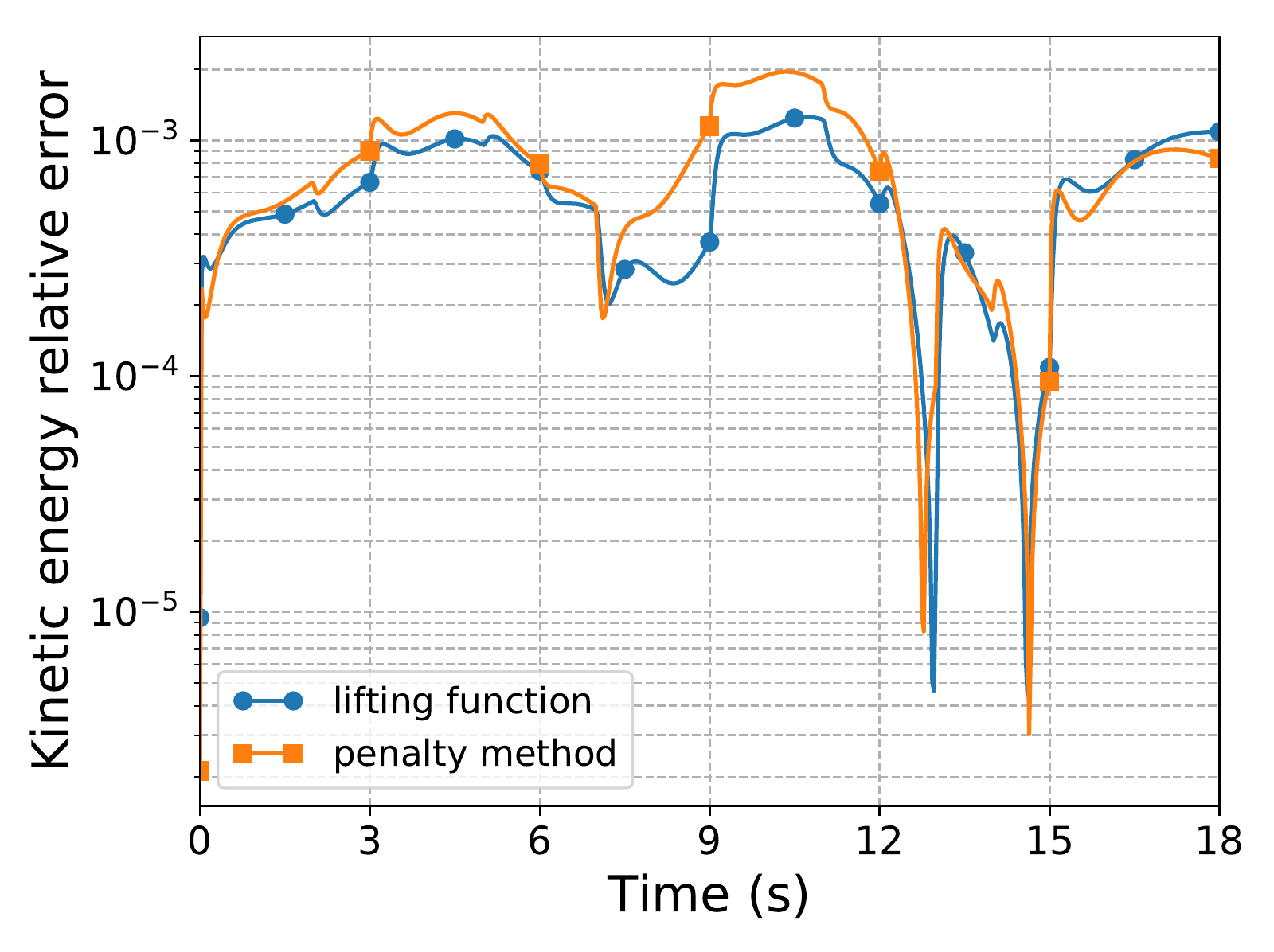}
	\caption{Kinetic energy relative $L^2$-error for the ROM with lifting function and with penalty method for the Y-junction flow problem.}
	\label{fig:Y_KE_eps}
\end{figure}

The relative error is more or less the same for both boundary control methods, as also was observed previously for the lid driven cavity test case, except around 9 s of simulation time. Then the difference in relative error for pressure between the two methods is the largest; the penalty method is about 2 $\cdot$ \num{e-1} larger than the error obtained with the lifting function method. However, on the long term the penalty method performs slightly better. This can also be concluded by having a look at the kinetic energy relative error in Figure~\ref{fig:Y_KE_eps}. Other than that, the relative velocity error is of the order \num{e-2} and for pressure \num{e-1}.
A possible source for the larger pressure error is that the PIMPLE algorithm, consisting of predictor and correction steps for pressure and velocity, is used at full order level, while the coupled (pressure-velocity) system at reduced order level is solved with Newton's iterative method. This is causing a discrepancy between the full order and reduced order model formulation. Nevertheless, the difference between the minimum and maximum relative error for both variables is about one order. 

\begin{figure}[h!]
	\centering
	\begin{subfigure}{.19\linewidth}
		\centering
		\includegraphics[width=1.\linewidth]{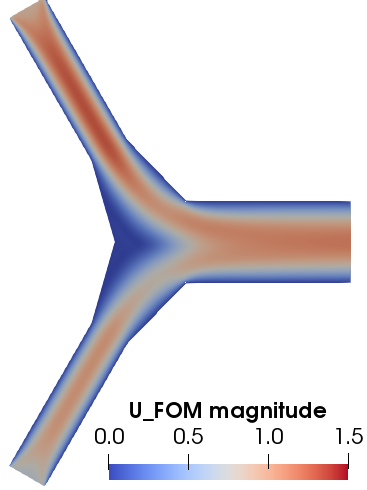}
	\end{subfigure}%
	\begin{subfigure}{.19\textwidth}
		\centering
		\includegraphics[width=1.\linewidth]{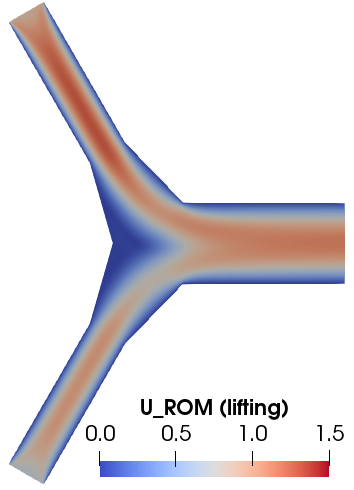}
	\end{subfigure}
	\begin{subfigure}{.19\textwidth}
		\centering
		\includegraphics[width=1.\linewidth]{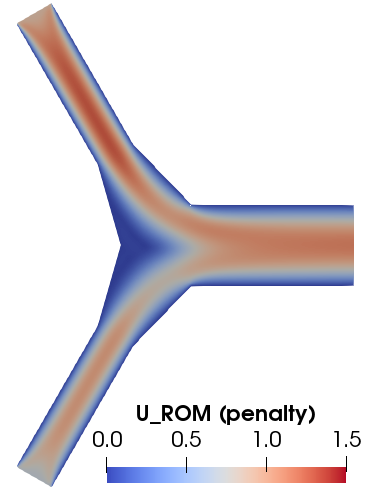}
	\end{subfigure}
	\begin{subfigure}{.19\textwidth}
		\centering
		\includegraphics[width=1.\linewidth]{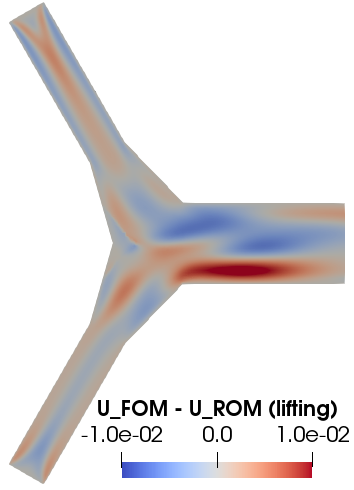}
	\end{subfigure}
	\begin{subfigure}{.19\textwidth}
		\centering
		\includegraphics[width=1.\linewidth]{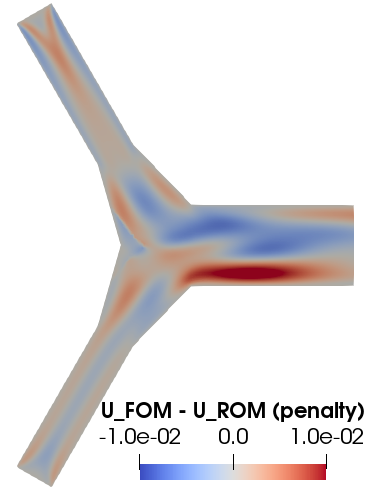}
	\end{subfigure}
	
	\begin{subfigure}{.19\linewidth}
		\centering
		\includegraphics[width=1.\linewidth]{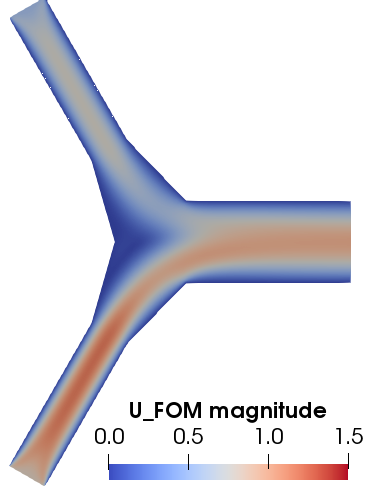}
	\end{subfigure}%
	\begin{subfigure}{.19\textwidth}
		\centering
		\includegraphics[width=1.\linewidth]{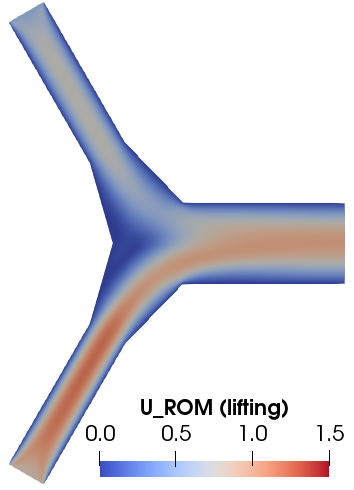}
	\end{subfigure}
	\begin{subfigure}{.19\textwidth}
		\centering
		\includegraphics[width=1.\linewidth]{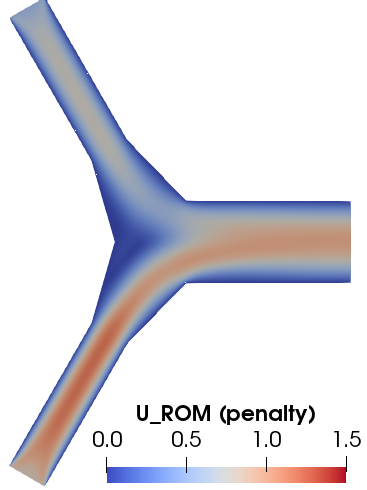}
	\end{subfigure}
	\begin{subfigure}{.19\textwidth}
		\centering
		\includegraphics[width=1.\linewidth]{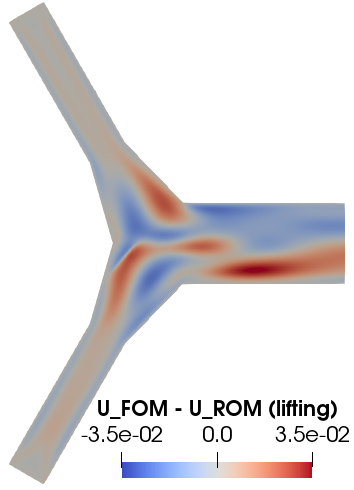}
	\end{subfigure}
	\begin{subfigure}{.19\textwidth}
		\centering
		\includegraphics[width=1.\linewidth]{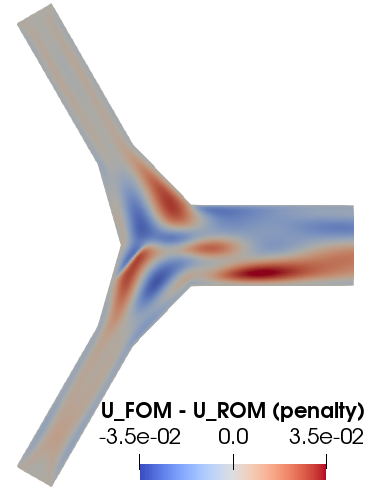}
	\end{subfigure}
	
	\begin{subfigure}{.19\linewidth}
		\centering
		\includegraphics[width=1.\linewidth]{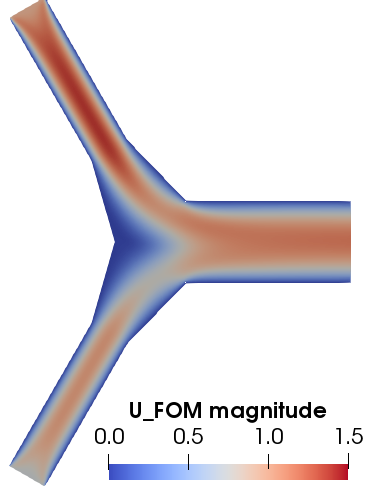}
	\end{subfigure}%
	\begin{subfigure}{.19\textwidth}
		\centering
		\includegraphics[width=1.\linewidth]{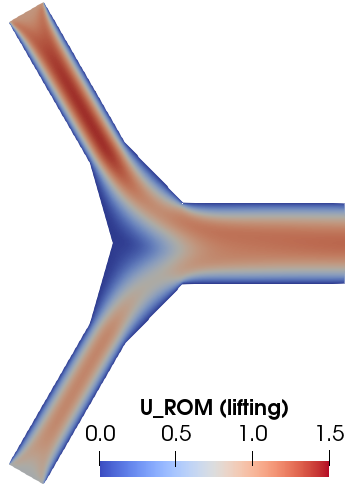}
	\end{subfigure}
	\begin{subfigure}{.19\textwidth}
		\centering
		\includegraphics[width=1.\linewidth]{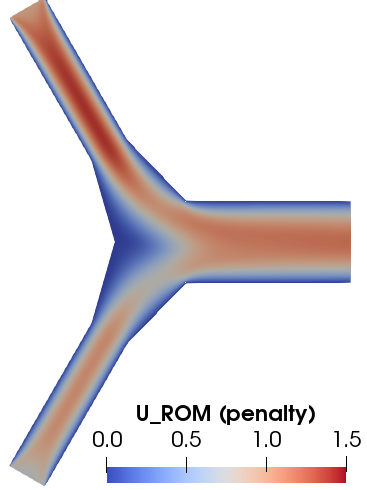}
	\end{subfigure}
	\begin{subfigure}{.19\textwidth}
		\centering
		\includegraphics[width=1.\linewidth]{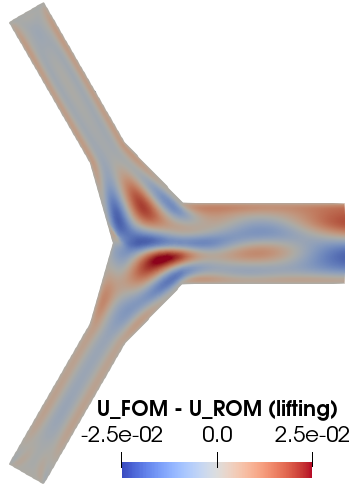}
	\end{subfigure}
	\begin{subfigure}{.19\textwidth}
		\centering
		\includegraphics[width=1.\linewidth]{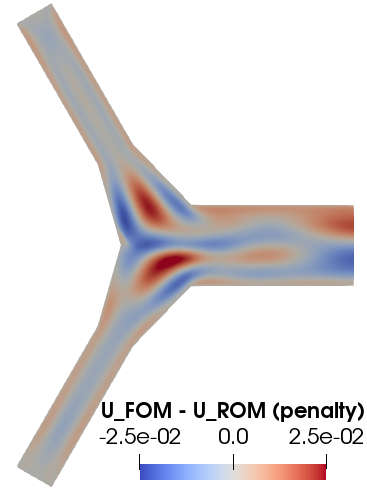}
	\end{subfigure}
	
	\caption{Comparison of the full order velocity magnitude fields (1st column), the ROM fields obtained with the lifting function method (2nd column) and penalty method (4th column) and the difference between the FOM and ROM fields obtained with the lifting function method (3rd column) and penalty method (5th column) at $t$ = 3, 9 and 18 s (from top to bottom) for the Y-junction flow problem. }
	\label{fig:Y_U}
\end{figure}

\begin{figure}[h!]
	\centering
	\begin{subfigure}{.19\linewidth}
		\centering
		\includegraphics[width=1.\linewidth]{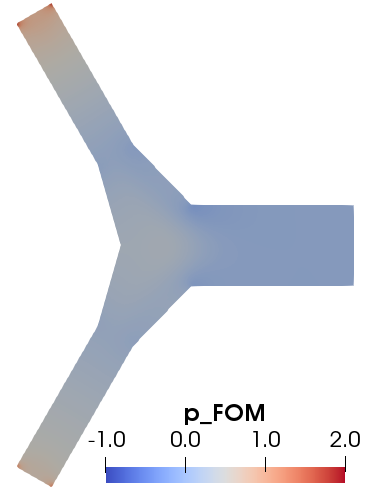}
	\end{subfigure}%
	\begin{subfigure}{.19\textwidth}
		\centering
		\includegraphics[width=1.\linewidth]{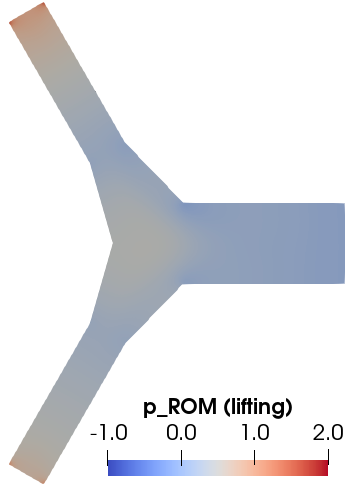}
	\end{subfigure}
	\begin{subfigure}{.19\textwidth}
		\centering
		\includegraphics[width=1.\linewidth]{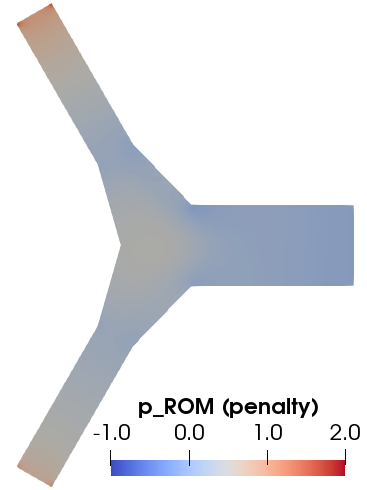}
	\end{subfigure}
	\begin{subfigure}{.19\textwidth}
		\centering
		\includegraphics[width=1.\linewidth]{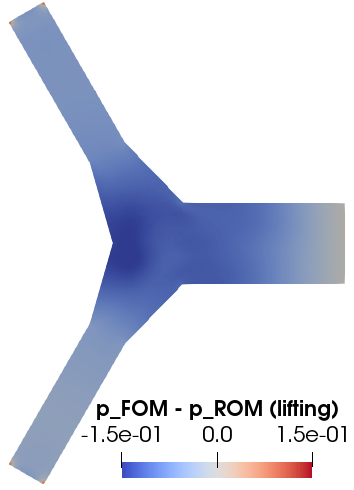}
	\end{subfigure}
	\begin{subfigure}{.19\textwidth}
		\centering
		\includegraphics[width=1.\linewidth]{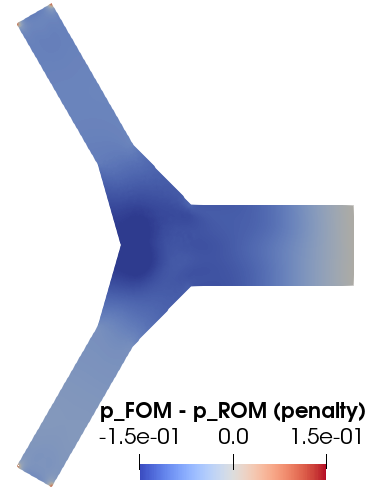}
	\end{subfigure}
	
	\begin{subfigure}{.19\linewidth}
		\centering
		\includegraphics[width=1.\linewidth]{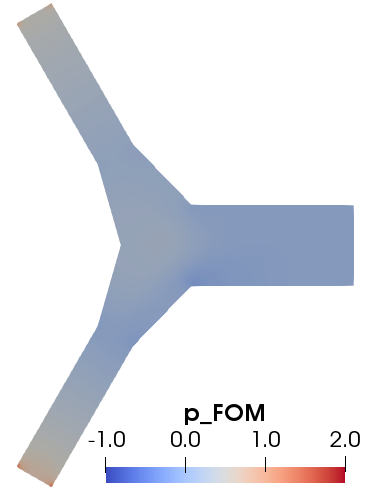}
	\end{subfigure}%
	\begin{subfigure}{.19\textwidth}
		\centering
		\includegraphics[width=1.\linewidth]{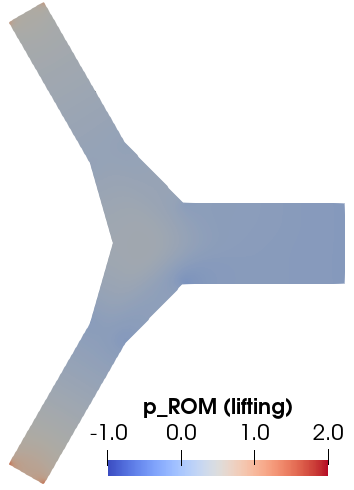}
	\end{subfigure}
	\begin{subfigure}{.19\textwidth}
		\centering
		\includegraphics[width=1.\linewidth]{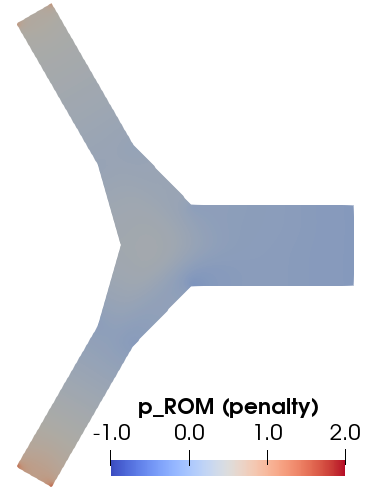}
	\end{subfigure}
	\begin{subfigure}{.19\textwidth}
		\centering
		\includegraphics[width=1.\linewidth]{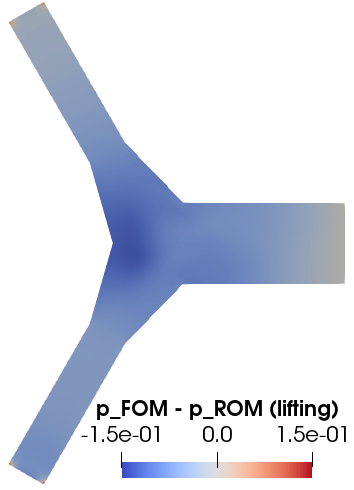}
	\end{subfigure}
	\begin{subfigure}{.19\textwidth}
		\centering
		\includegraphics[width=1.\linewidth]{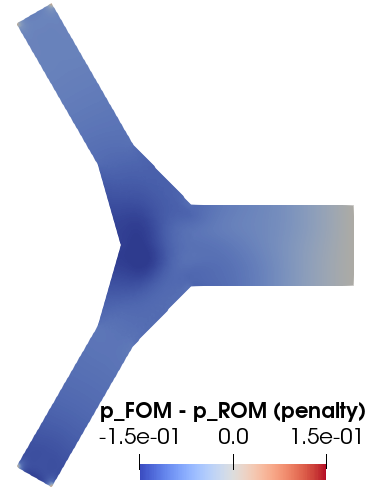}
	\end{subfigure}
	
	\begin{subfigure}{.19\linewidth}
		\centering
		\includegraphics[width=1.\linewidth]{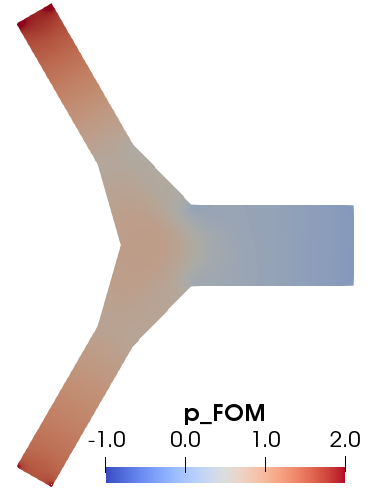}
	\end{subfigure}%
	\begin{subfigure}{.19\textwidth}
		\centering
		\includegraphics[width=1.\linewidth]{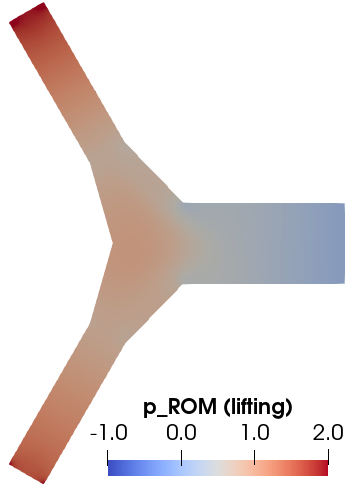}
	\end{subfigure}
	\begin{subfigure}{.19\textwidth}
		\centering
		\includegraphics[width=1.\linewidth]{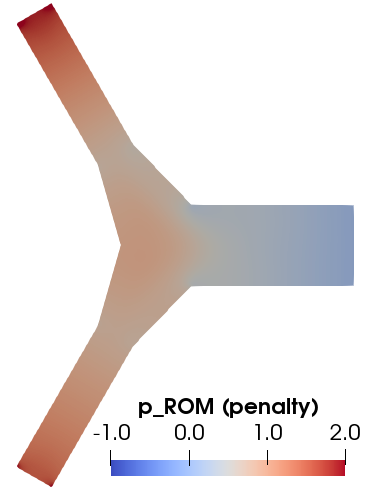}
	\end{subfigure}
	\begin{subfigure}{.19\textwidth}
		\centering
		\includegraphics[width=1.\linewidth]{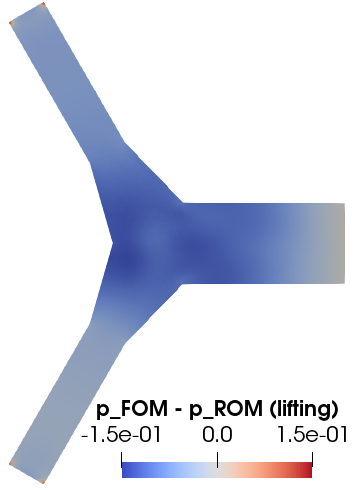}
	\end{subfigure}
	\begin{subfigure}{.19\textwidth}
		\centering
		\includegraphics[width=1.\linewidth]{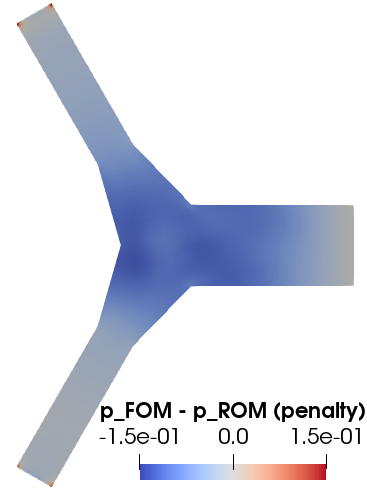}
	\end{subfigure}
	
	\caption{Comparison of the full order pressure fields (1st column), the ROM fields obtained with the lifting function method (2nd column) and penalty method (4th column) and the difference between the FOM and ROM fields obtained with the lifting function method (3rd column) and penalty method (5th column) at $t$ = 3, 9 and 18 s (from top to bottom) for the Y-junction flow problem. }
	\label{fig:Y_p}
\end{figure}

Furthermore, the absolute error between the FOM and the ROMs is shown in Figures~\ref{fig:Y_U} and~\ref{fig:Y_p} for velocity magnitude and pressure, respectively. For velocity the absolute error between the FOM and the ROM is of the order \num{e-2} for all plotted simulation times and the absolute error for pressure is of the order \num{e-1}. For pressure, the error is indeed larger in the case of the penalty method compared to the lifting function method at 9 s of simulation time, as previously observed in Figure~\ref{fig:Y_L2_error}, but in general, the error distribution, for both the velocity and pressure fields, is similarly distributed over the domain, and thus the methods are performing the same.

\clearpage
Finally, the computational times for performing the full order simulation (Eq.~\ref{eq:FOM_mat}), calculating the POD modes (Eqs.~\ref{eq:approx}-~\ref{eq:POD}), the reduced matrices (Eqs.~\ref{eq:ROM_matrices},~\ref{eq:C_matrix} and~\ref{eq:ROM_matrices2}) and performing the simulation at reduced level (Eq.~\ref{eq:ROM} (lifting function method) or Eq.~\ref{eq:pen_ROM_min} (penalty method) \& Eq.~\ref{eq:ROM2}) are listed in Table~\ref{tab:Y_times}. Calculating the reduced matrices and the ROM solutions takes more time in the case of the lifting function method as the reduced basis space consists of four additional modes, namely the normalized lifting functions, compared to the penalty method. Determining the penalty factor with the iterative method takes 1.4 s. The speedup ratio between the ROM and the FOM is about 13 times for the lifting method and 24 times for the iterative penalty method. 

\begin{table}[h!]
	\caption{Computational time (clock time) for the FOM simulation, POD modes, calculating reduced matrices offline (Matrices), determining penalty factor with iterative method (Penalty factor) and ROM simulation.}
	\centering
	\begin{tabular}{lllllll}
		\hline
		\multicolumn{1}{c}{Method}  &FOM  &POD  &Matrices &Penalty factor  &ROM   \\ \hline
		Lifting                       &13 min.   & 7.6 s  & 9.2 s &-  &59 s    \\
		Penalty 	  & 13 min. & 7.9 s  & 4.7 s &1.4 s &33 s  \\ \hline
	\end{tabular}
	\label{tab:Y_times}
\end{table}

\section{DISCUSSION}\label{sec:discussion}
The results have shown that the lifting function method and penalty method perform equally and lead to similar results. However, they have their own advantages and drawbacks. A disadvantage of the penalty methods is that the penalty factor cannot be determined a priori~\cite{graham1999optimal1}. The implementation of an iterative solver to determine the penalty factor does however save time compared to performing numerical experimentation manually. On the other hand, even though a lifting function(s) can be determined beforehand, it may be hard to find a function that will lead to an accurate ROM and therefore extensive testing of ROMs for different functions can be needed. In this work, the lifting functions are obtained by solving a potential flow problem and are thus physics-based unlike the penalty factor, which is an arbitrary value. Moreover, this value needs to be chosen above a certain threshold to enforce the BCs in the ROM, but can lead to an inaccurate ROM solution if it is too high~\cite{Sirisup}. In that case, the penalty method fails for that specific problem. 

Finally, an advantage of the penalty method stated in literature~\cite{lorenzi2016pod} is that long-time integration and initial condition issues are less of a problem compared to a lifting function method. Here the ROMs have not been tested for long-term integration, so further research is needed in order to confirm this statement. However, as tested for the Y-junction test case, the ROM is accurate and does not exhibit instabilities even outside the time domain in which snapshots were collected. 

For both cases tested in this study, only one full order simulation has been performed for collecting the snapshots. However, in case the BCs of the Y-junction are not time-dependent, snapshots from at least two different offline solves are required for the penalty method. The reason for this is that the boundary conditions are a linear combination of snapshots and the boundary conditions can therefore only be scaled and not be set to any arbitrary value in case only snapshots from one full order simulation are used for the POD. When several sets of snapshots for different boundary values are required, one can optimize the POD procedure by using a nested POD approach~\cite{georgaka2018parametric}. 

It is important to note that the penalty factor is determined during the online phase and does not depend on high-fidelity data. Therefore, no modification are needed in the case of parametric problems that, for example, use the viscosity as the physical parameter.

In the case of the Y-junction test case, the penalty method can be used to adjust the direction of the inlet flow in the ROM. One penalty factor is enforcing the x-direction and another the y-direction. Nevertheless, new snapshots for different inlet angles are required as the current POD bases do not contain this information. For the lifting function method, it is often problematic to determine suitable lifting functions that are physical. Ideally, the lifting functions are orthogonal to each other as in the work of Hijazi et al.~\cite{HijaziStabileMolaRozza2020a} who studied a flow past an airfoil with parameterized angle of attack and inflow velocity. They used two lifting functions with orthogonal inflow conditions: $\boldsymbol{\zeta}_{c_1}$ = (0,1) and $\boldsymbol{\zeta}_{c_2}$ = (1,0) on ${\Gamma_i}$, respectively. These lifting functions are obtained by solving two linear potential flow problems. In that way, it is possible to adjust the direction of the flow at a inlet by scaling the associate lifting functions accordingly. However, specifying a purely tangential velocity at the inlets of the Y-junction would result in unphysical lifting functions. Thus, this approach is only suitable for a few problems and will not always lead to physical results. 

In the case of non-physical lifting functions, the ROM gets unstable or the ROM solutions are polluted with noise. This strongly depends on the chosen lifting functions.

Moreover, both methods can, in theory, also be used for controlling pressure boundary conditions, but this is not studied in this work. 

In this study, the exploitation of a pressure Poisson equation has been incorporated in the ROM as a stabilization method. Even though the ROMs are indeed stable, the relative error for pressure is about an order higher than for velocity. Alternatively, the supremizer enrichment of the velocity space technique could be used to stabilize the ROM, which may lead to more accurate pressure fields~\cite{Stabile2017CAF,kean2020error,ballarin2015supremizer}.

Furthermore, the ROMs can be improved by using a second order backward method for the time discretization of the ROM as the FOMs are treated using a second order backward differencing scheme.

Finally, for the Y-junction test case, the full order snapshots and the ROM solutions all have inlet velocities between an identical maximum and minimum value. The ROM could become less stable and accurate in case it is tested for values outside this range. Therefore it is recommended to collect snapshots for the same range as for which the ROM boundary needs to be controlled.

\section{CONCLUSIONS AND PERSPECTIVES}\label{sec:conclusion}
Two boundary control methods are tested: the lifting function method and the iterative penalty method for controlling the velocity boundary conditions of FV-based POD-Galerkin ROMs. The penalty method has been improved by using an iterative solver for the determination of the penalty factors, rather than using numerical experimentation. The factors are determined by the iterative solver in about a second for both test cases. The results of the reconstructed velocity and pressure fields show that both methods are performing equally. Moreover, the reduced order model of which the boundary conditions are controlled with the iterative penalty method is about two times faster compared to the lifting function method for the Y-junction flow case.   

A pressure Poisson equation approach is applied for the reconstruction of the pressure field and to stabilize the ROM. For time-dependent boundary problems, an additional term is added to the ROM formulation.

Finally, a speedup factor, the ratio between the FOM and ROM simulation time, of 308 is obtained with the iterative penalty method and of 270 with the lifting function method for the lid driven cavity test case. The speedup factors are 24 and 13, respectively, for the Y-junction test case.

For further development, the model will be extended for turbulent flows, which will be essential to simulate industrial flow problems. Furthermore, the control of the pressure boundary conditions needs to be investigated, which may be required when coupling 3D CFD problems with 1D system codes. Also, the accuracy of the reconstructed pressure fields can be improved by using a supremizer enrichment approach rather than solving the Pressure Poisson Equation. The effect of supremizer enrichment of the velocity space on the boundary control methods will have to be investigated.

\section*{Acknowledgments} 
This work has been partially supported by the ENEN+ project that has received funding from the Euratom research and training Work Programme 2016 - 2017 - 1 \#755576. In addition we acknowledge the support provided by the European Research Council Executive Agency by the Consolidator Grant project AROMA-CFD ``Advanced Reduced Order Methods with Applications in Computational Fluid Dynamics'' - GA 681447, H2020-ERC CoG 2015 AROMA-CFD and INdAM-GNCS projects.

\bibliographystyle{ieeetr}
\bibliography{MyBib}

\end{document}